\documentclass[final,leqno,onefignum,onetabnum]{siamltex}
\usepackage{fullpage}
\usepackage{comment}
\usepackage{upgreek}
\usepackage[mathscr]{eucal}
\usepackage{color}
\usepackage[usenames,dvipsnames,svgnames,table]{xcolor}
\usepackage{amsmath,amssymb,amsopn,mathtools}
\usepackage{graphicx}
\usepackage{hyperref}       
\hypersetup{
	colorlinks=true,   
	citecolor=	PineGreen,     
	filecolor=PineGreen,     
	linkcolor=PineGreen,    
	urlcolor= PineGreen}      

\usepackage{relsize}
\numberwithin{equation}{section}
\usepackage{enumitem}
\newtheorem{remark}[theorem]{Remark}

   \allowdisplaybreaks
   \DeclareMathOperator*{\argmin}{argmin}

\newcommand{\tu}{\textup}

\newcommand{\dd}{\displaystyle}
\newcommand{\bfa}[1]{\boldsymbol{#1}} 			%

\usepackage{color}
\usepackage{adjustbox}
\usepackage{diagbox}
\definecolor{black}{rgb}{0,0,0}

\definecolor{red}{rgb}{1,0,0}

\definecolor{blue}{rgb}{0,0,1}

\usepackage{subcaption}
\renewcommand{\theequation}{\arabic{section}.\arabic{equation}}
\numberwithin{equation}{section}

\renewcommand{\div}{\mathop{\rm div}\nolimits}

\newcommand{\todo}[1]{{\color{red}{#1}}}
\newcommand{\jrp}[1]{{\color{orange}{#1}}}
\newcommand{\tony}[1]{{\color{blue}{#1}}}

\newcommand{\dx}{ \mathrm{d}x}

\newcommand{\dt}{ \mathrm{d}t}

\newcommand{\beq}{\begin{equation}}
\newcommand{\eeq}{\end{equation}}
\newcommand{\beqq}{\begin{equation*}}
\newcommand{\eeqq}{\end{equation*}}
\newcommand{\beqas}{\begin{eqnarray*}}
\newcommand{\eeqas}{\end{eqnarray*}}
\newcommand{\bsp}{\begin{split}}
\newcommand{\eesp}{\end{split}}

\usepackage[english]{babel}

\makeatletter
\def\thanks#1{\protected@xdef\@thanks{\@thanks
		\protect\footnotetext{#1}}}
\makeatother

\begin{document}
\title{Constraint Energy Minimizing Generalized Multiscale Finite Element Method for multi-continuum Richards equations}
\author{
Tina Mai$^{\MakeLowercase{a,b}}$\thanks{\textit{Tina Mai}; $^a$Institute of Research and Development, Duy Tan University, Da Nang, 550000, Vietnam; $^b$Faculty of Natural Sciences, Duy Tan University, Da Nang, 550000, Vietnam; \texttt{maitina@duytan.edu.vn}
}
\and
Siu Wun Cheung$^{\MakeLowercase{c}}$\thanks{\textit{Siu Wun Cheung}; $^c$Center for Applied Scientific Computing,
	Lawrence Livermore National Laboratory, Livermore, CA 94550, USA; \texttt{cheung26@llnl.gov}
}
\and
Jun Sur Richard Park$^{\MakeLowercase{d,*}}$\thanks{$^{*}$Corresponding author: \textit{Jun Sur Richard Park}; $^d$Department of Mathematics, The University of Iowa, Iowa City, IA, USA; \texttt{junsur-park@uiowa.edu}
}
}
\maketitle

\maketitle

\begin{abstract}
In fluid flow simulation, the multi-continuum model is a useful strategy.  When the heterogeneity and contrast of coefficients are high, the system becomes multiscale, and some kinds of reduced order methods are demanded.  Combining these techniques with nonlinearity, we will consider in this paper a dual-continuum model which is generalized as
a multi-continuum model for a coupled system of nonlinear Richards equations as unsaturated flows, in complex heterogeneous fractured porous media; and we will solve it by a novel multiscale approach utilizing
the constraint energy minimizing generalized multiscale finite element method (CEM-GMsFEM). In particular, such a nonlinear system will
be discretized in time and then linearized by Picard iteration (whose global convergence is proved theoretically). Subsequently,
we tackle the resulting linearized equations by the CEM-GMsFEM and obtain proper offline multiscale basis functions to span the multiscale space (which contains the pressure solution).  More specifically, we first introduce two new sources of samples, and the GMsFEM is used over each coarse block to build local auxiliary multiscale basis functions via solving local spectral problems, that are crucial for detecting high-contrast channels. Second, per oversampled coarse region, local multiscale basis functions are created through the CEM as constrainedly minimizing an energy functional.  
Various numerical tests for our approach reveal that the error converges with the coarse-grid size alone and that only few oversampling layers as well as basis functions are needed.
%

\end{abstract}



\section{Introduction}
For any soil sample, the amount of water retained within the gaps between unsaturated soil particles is known as soil moisture.  Even being a small portion in many parts of the water cycle, soil moisture is
crucial to various procedures of hydrology, biology, and biogeochemistry.  For example, soil moisture is a key variable to farming, environmental management, groundwater storage, geotechnics, energy balances, meteorological forecast, and earth system dynamics, etc.  Richards equation \cite{r0,richarde1,richarde2,richarde3,richardsreview}, which features the seepage of water into some porous material 
filled with water and air \cite{ry1}, is used as an unsaturated flow to quantitatively model the associated processes.  Evaporation and precipitation, which are tightly coupled in nonlinear ways, affect moisture near the soil surface the most, prompting us to explore a coupled system of nonlinear Richards equations.

Also, in our considering porous media, there can exist complex heterogeneous rock properties, faults, intricate fracture geometries, multi-continuum background with mass transfer, high contrast, and numerous scales, among other aspects.  Especially, the material characteristics of fractures can differ significantly from those of the surrounding media, which can also comprise extremely heterogeneous and large-contrast regions as well as high permeability.  
These obstacles lead to the fact that they can have a considerable impact on nonlinear fluid flow processes and solutions comprise multiple scales, 
making traditional numerical simulations much more difficult because extra computing power is needed.

The purpose of this study is to build and examine some reduced models for these types of issues.  In the standard upscaling methods through homogenization, the computational domain is first partitioned into coarse-scale blocks, where scales are not necessarily resolved, then effective material property for each coarse block are calculated employing the fine-scale solutions of some local problem \cite{dur91,weh02}.
However, it is well understood that one effective coefficient per coarse patch is insufficient to represent all features of the solutions, particularly in the regions holding important modes, fractures, high-contrast heterogeneities, and interaction of continua.

To resolve this disadvantage, we utilize on coarse grid the multi-continuum strategies \cite{baren,arbogast20,warren1963behavior,kazemi1976numerical,wu1988multiple,pruess1982fluid}, where a number of effective medium properties are built.  Physically, each continuum is treated as a system (throughout the entire domain) so that the flow between them can be easily characterized.  Different continua are adjacent in the fine grid. They coexist via mean characteristics \cite{baren} at every location of the considering region on the coarse grid, and interactions appear among them.  Mathematically, we represent on each coarse block a system of equations, each of which corresponds to one of the fine grid's multicontinua.  In this paper, using dual-continuum model for the unsaturated flows, we construct distinct Richards equations for the flow in natural fractures and the flow within matrix (background), and some specific interaction terms are coupled to such equations, as in \cite{baren, douglas1990dual, warren1963behavior}.  This purpose is achieved by assuming that each continuum is connected to the other (even if it is not topologically connected, across the kind of coupling and the entire domain), provided that it possesses global effects solely. 

To illustrate our multi-continuum strategy, we now look at dual-continuum background in further detail. Barenblatt \cite{baren} developed the first dual-porosity model for flow simulation in fissured rock.  The proposed two continua in that work are for characterizing low and high porosity continua, namely, a system of natural fractures (so-called small-scale connected, highly developed, or well-developed fractures) and a matrix, both of which are used in our paper. On the basis of \cite{baren}, there was also an early work using homogenization on dual continua \cite{arbogast20}.
Intraflow and interflow transfers are together considered per continuum. Essentially, the dual-continuum background can take arbitrary shape and fit any of the above approaches.


Dual-continuum models are also utilized to represent a variety of scientific and engineering applications, such as complicated processes in shale reservoirs \cite{akkutlu2012multiscale,akkutlu2017,organic17},
where those models are employed to depict a complex interplay of organic and inorganic matter. Also, dual-continuum models can characterize flow through vugs 
and the rest media in vuggy carbonate reservoirs' simulations \cite{wu2006multiple, wu2011multiple, Tony11}.



The classical direct approach to tackle multi-continuum models with fractures is local fine-grid simulation, in a few simple steps \cite{karimi2016general}.
First, a fine grid is built locally to represent the shapes of fractures and heterogeneities of background. Second, the flow equations are discretized on that fine grid, and a global solution is obtained from the collection of local solutions.  This technique can be implemented using well-known frameworks, such as the Finite Element Method (FEM) \cite{Baca84} and the Finite Volume Method (FVM) \cite{bogdanov2003two, granet2001two, karimi2004efficient, monteagudo2004control, noetinger2015quasi, reichenberger2006mixed}.  Within the confines of the finite-element framework, there are considerations of the ordinary Galerkin formulation in \cite{Baca84, juanes2002general, karimi2003numerical, kim2000finite}, the mixed finite element approach in \cite{erhel2009flow, hoteit2008efficient, ma2006mixed, martin2005modeling}, the hierarchical FEM
in \cite{rh1}, and the discontinuous Galerkin method in \cite{eikemo2009discontinuous, hoteit2005multicomponent}.  A hybrid strategy has also been studied \cite{geiger2009black, matthai2007finite, nick2011comparison}, which combines the FVM for the transport equation with the FEM for the pressure equation.  However, even with the aid of supercomputers and parallel computing, direct fine-grid simulation of multiple-scale problems is difficult and expensive, leading to the need of some multiscale methods.  

The inspiration for the novel multiscale approach we develop in this article is the generalized multiscale finite element method (GMsFEM) \cite{G1,chung2016adaptive,chung2015generalizedwave,cho2017generalized}, which may be thought of as a generalization of the multiscale finite element method (MsFEM) \cite{Ms,Msnon}.  We will build coarse-grid multiscale basis functions that can couple multiple continua together with high-contrast channels, to achieve the small-scale impact on the large scales without the urge to solve for all minor intricacies.  The GMsFEM's primary idea is to employ local spectral decomposition in some appropriate snapshot spaces to find local dominant modes.  The resulting dominant eigenfunctions can transmit local to global properties using coarse-grid multiscale basis functions.  These concepts are crucial for recognizing the effects of high-contrast regions as well as channels, which must be described separately by distinct basis functions.  
For instance, if we have $n$ different connected fracture networks within a coarse block, then there are $n$ very small eigenvalues, and the related dominant eigenvectors will reflect these connected fracture networks and can be regarded as lowered degrees of freedom that depict these fracture effects.
In this way, the GMsFEM and multi-continuum techniques possess certain commonalities (see \cite{mcontinua17}, for instance).  A variety of domain decomposition methods \cite{pre2,kim2015bddc,kim2016bddc} have leveraged the idea of building local basis functions utilizing spectral problem. 
Recently, the GMsFEM has been successfully applied to a variety of problems \cite{organic17,Tony11,mcontinua17, Spiridonov2019, ericmultiporoelastic19a,akkutlu2017,akkutlu-poro-25,cosserat2022,yanbo1,yanbo2,yanbo3, yanbo4}.

In our previous work \cite{rtt21}, the GMsFEM was utilized to solve dual-continuum Richards equations within complex heterogeneous fractured porous media.  Nevertheless, it is not straightforward to construct a multiscale approach whose convergence is only determined by the coarse-grid size and is unaffected by scales or contrast.  Several ways are discussed in the literature to generate multiscale algorithms with mesh-dependent convergence \cite{owhadi2014polyharmonic, maalqvist2014localization, owhadi2017multigrid,hou2017sparse,cem1,chung2018mixed}. The use of local spectral problems to capture the impact of high-contrast channels is motivated by the GMsFEM's theory, and this principle is also applied to mesh-dependent convergence \cite{hou2017sparse,cem1,chung2018mixed}.

\newpage

In that spirit, a new multiscale technique for a linear dual-continuum model was built and investigated in \cite{cheung2018constraint}, following the constraint energy minimizing generalized multiscale finite element method (CEM-GMsFEM) \cite{cem1,chung2018mixed,icem}.  Such technique relies on the coarse-grid size alone for convergence. 
In our paper, for the case of multi-continuum nonlinear Richards equations, after temporal discretization, at each time step until the halting time, we employ linearization in Picard iteration (with a desired ending indicator), and the CEM-GMsFEM \cite{cheung2018constraint} is applied to find the solution of these linearized equations in two-dimensional porous fractured heterogeneous media.  In Appendix \ref{cp}, we shall additionally show that the Picard iteration process converges globally. 

The CEM-GMsFEM in this paper is made up of two components: the creation of local basis functions for each coarse element (via employing the GMsFEM to build the auxiliary multiscale basis functions) and then per oversampled coarse domain (by using the CEM to achieve the collection of multiscale basis functions with locally minimal energy).  More specifically, after introducing two new sources of samples, we will first establish local auxiliary multiscale basis functions over each coarse block through the GMsFEM.  The number of such functions is the same as the number of high-contrast channels.  These functions are dominant eigenfunctions 
(corresponding to the smallest eigenvalues of local spectral problems) and can be considered as the reduced degrees of freedom needed to model channelized effects. Such eigenfunctions are also important in the establishment of localized basis functions, as we point out.  Second, multiscale basis functions form the other important component.  For each oversampled coarse region, these functions are created by the CEM as minimizing an energy functional, which is constrained in a manner that its minimizer satisfies a set of orthogonality requirements with respect to the auxiliary functions. 
It is clear from the numerical results that the error converges with the coarse-grid size only, so our multiscale approach derived by a Galerkin formulation has the mesh-dependent convergence property.

Adaptivity can also be performed in the CEM-GMsFEM, as shown in \cite{chungres1,chung2015residual,cem2f}.  Some recent applications of the CEM-GMsFEM can be found in \cite{ye22cem,ec22cem,cemDGelastic}.  With the CEM-GMsFEM's development, there have been important studies on the non-local multi-continuum (NLMC) method  \cite{nlmcdual,vasilyeva2019nonlocal,chung2018non,cem28,nlmc14}. The NLMC method's main principle is similar to the CEM-GMsFEM, with an exception that the multiscale basis functions are changed to reflect the average of solutions such that the degrees of freedom possess physical meanings.  In multi-continuum fractured media, these techniques are also effective in tackling high-contrast and multiscale components.

The paper is structured in the following way.  
In Section~\ref{sec:model}, the multi-continuum model for a coupled system of nonlinear Richards equations is introduced, within fractured heterogeneous porous media.  We present in Section \ref{pre} the fine-scale discretization and Picard iteration for linearization of such system.  In Section~\ref{sec:method}, our novel multiscale approach will be provided to solve this linearized system, using a new idea of two sample sources for pressure in constructing multiscale spaces, following the constraint energy minimizing generalized multiscale finite element method (CEM-GMsFEM).
In Section~\ref{sec:numerical}, various numerical examples will be shown to expose the approach's mesh-dependent convergence. The paper is concluded in Section~\ref{sec:conclusions}.  We give a proof for the Picard linearization's global convergence in Appendix \ref{cp}.

\section{Multi-continuum Richards equations}\label{sec:model}
Let $\Omega$ be a bounded, simply connected, open, Lipschitz, convex computational domain in $\mathbb{R}^d$.
The case $d=2$ is considered to ease our discussion throughout this paper, however the method can be easily generalized to $d=3\,.$
The subscripts $i$ and $l$ stand for indices of continua, 
where $N$ denotes the number of continua.  The symbols $\dfrac{\partial}{\partial t}$ and $\nabla$ respectively represent the temporal derivative and spatial gradient.  Other notation is as in \cite{mcl,rtt21}.  
Vector fields and matrix fields over $\Omega$ are denoted by bold letters (e.g., $\bfa{v}$ and $\bfa{T}$) while functions are represented by italic capitals (e.g., $f$).  Over $\Omega\,,$ the spaces of functions, vector fields, and matrix fields are respectively expressed by italic capitals (e.g., $L^2(\Omega)$), 
boldface Roman capitals (e.g., $\bfa{V}$), 
and special Roman capitals (e.g., $\mathbb{S}$). 


At first glance, a coupled system of dual-continuum nonlinear Richards equations
\cite{rtt21,Spiridonov2019} has the form
\beq
\label{eq:original1}
\bsp
\frac{\partial p_1(t,\bfa{x})}{\partial t} - \div [\varkappa_1(\bfa{x},p_1(t,\bfa{x}))\nabla p_1(t,\bfa{x})]+ Q_{12}(\bfa{x},p_1(t,\bfa{x})&, p_2(t,\bfa{x}))(p_1(t,\bfa{x})-p_2 (t,\bfa{x}))\\
& = f_1(t,\bfa{x}) \  \textrm{in} \ (0,T] \times \Omega\,,\\
\frac{\partial p_2(t,\bfa{x})}{\partial t} - \div [\varkappa_2(\bfa{x},p_2(t,\bfa{x}))\nabla p_2(t,\bfa{x})] + Q_{21}(\bfa{x},p_2(t,\bfa{x})&, p_1(t,\bfa{x}))(p_2(t,\bfa{x})-p_1(t,\bfa{x}))\\
& = f_2(t, \bfa{x}) \ \textrm{in} \ (0,T] \times \Omega\,.
\end{split}
\eeq
In this paper, we consider a general multi-continuum model of such system as follows \cite{Spiridonov2019,nlnlmc31}: for each continuum $i = 1,\ldots,N\,,$ 
%
\beq
\label{eq:original}
\bsp
\frac{\partial p_i}{\partial t} - \div (\varkappa_i(p_i)\nabla p_i)+  \sum_{l=1}^N Q_{il}(\bfa{x},p_i, p_l)(p_i-p_{l}) = f_i(t,\bfa{x}) \  \textrm{in} \ (0,T] \times \Omega\,,
\end{split}
\eeq
where $T>0$ is the final time.  This system is prescribed with the initial condition $p_i(0,\bfa{x})= p_{i,0}$ in $\Omega$ and the Dirichlet boundary condition $p_i(t,\bfa{x})=0$ on $(0,T] \times \partial \Omega\,.$  Basic notation can be found in \cite{rtt21}.
Here, $p_i:=p_i(t,\bfa{x})$ stands for the pressure head, 
$\varkappa_i(p_i):=\varkappa_i(\bfa{x},p_i)$ denotes the unsaturated hydraulic conductivity, $f_i$ represents the source or sink function for the $i$th continuum, and the term $Q_{il}(\bfa{x},p_i, p_l)(p_i - p_l)$ describes mass transfer  
of the liquid which flows from the $i$th continuum into the $l$th continuum per unit of media volume as well as per unit of time \cite{baren}, where we denote $Q_{il} = Q_{il}(p_i,p_l):=Q_{il}(\bfa{x},p_i,p_l)\,.$  When this mass exchange term $Q_{il}(\bfa{x},p_i, p_l)(p_i - p_l)$ vanishes, the system \eqref{eq:original} becomes a single-continuum equation (see Section \ref{sec:numerical}).

The $L^2$ inner product is represented by $(\cdot,\cdot)\,$, and the Sobolev space 
$V: = H_0^1(\Omega) = W_0^{1,2}(\Omega)\,$ 
is equipped with the norm $\| \cdot \|_{V}$: 
\[\|v\|_{V} = \left(\|v\|^2_{L^2(\Omega)} + 
\|\nabla v \|^2_{\bfa{L}^2(\Omega)}\right)^{\frac{1}{2}}\,.\]
Here, $\| \nabla v \|_{\bfa{L}^2(\Omega)}:= \| | \nabla v | \|_{L^2(\Omega)}\,,$ where 
$| \nabla v|$ indicates the Euclidean norm of the $d$-component vector-valued function 
$ \nabla v\,.$  We also denote $\bfa{V} = V^N = [H_0^1(\Omega)]^N\,.$  With $\bfa{v} = (v_1, \cdots, v_N) \in \bfa{V}\,,$ $\| \nabla \bfa{v}\|_{\mathbb{L}^2(\Omega)}:= \| | \nabla \bfa{v}| \|_{L^2(\Omega)}\,,$ where 
$| \nabla \bfa{v}|$ stands for the Frobenius norm of the $N \times d$ matrix $\nabla \bfa{v}\,.$ 

The hydraulic conductivity together with its spatial gradient as well as the mass transfer coefficient are assumed to be uniformly bounded, that is, positive constants $\underline{\varkappa}, \overline{\varkappa}$ and $\underline{\beta},\overline{\beta}$ exist so that the following inequalities are satisfied: 
\begin{align}
\label{Coercivity}
\begin{split}
\underline{\varkappa} \leq \varkappa_i(\bfa{x}, p_i), \ | \nabla \varkappa_i(\bfa{x},p_i)| \leq \overline{\varkappa}\,,\\
\underline{\beta} \leq Q_{il}(\bfa{x},p_i,p_l) \leq \overline{\beta}\,.
\end{split}
\end{align}

Without loss of generality, the initial condition is assumed to be 
\begin{equation}\label{ini}
\bfa{p}_0 = \bfa{p}(0,\bfa{x})= (p_{1,0},\ldots,p_{N,0}) \in \bfa{V}\,.
\end{equation}
Given $\bfa{u} = (u_1,\ldots,u_N) \in \bfa{V}$, with $i=1,2,\ldots,N$, we define the following bilinear forms: 
for all $\bfa{p}=(p_1,\cdots,p_N)\,, \ \bfa{v} = (v_1,\cdots,v_N)  \in \bfa{V}\,,$
\begin{align}
a_i(p_i,v_i;u_i)&=\int_{\Omega} \varkappa_i(u_i)\nabla p_i \cdot \nabla v_i \, \dx\,, \label{ai}\\
q_i(\bfa{p},v_i;\bfa{u})&= \sum_{l=1}^N \int_{\Omega} Q_{il}(u_i,u_l)(p_i-p_{l}) v_i \, \dx\,. \label{qi}
\end{align}
The variational form of \eqref{eq:original} reads:  
find $\bfa{p}=(p_1,\ldots,p_N) \in \bfa{V}$ such that with $i=1,\ldots,N$, 
\begin{align}\label{r1e}
\left(\frac{\partial{p_i}}{\partial t} , v_i \right) + a_i(p_i,v_i;p_i) + q_i(\bfa{p},v_i;\bfa{p}) = (f_i,v_i)\,,
\end{align}
for any $\bfa{v} = (v_1,\cdots,v_N) \in \bfa{V}\,,$ with a.e.\ $t \in (0,T]\,,$ and $f_i(t,\cdot) \in L^2(\Omega)\,.$  The initial condition is given in \eqref{ini}.
\section{Fine-scale discretization and Picard iteration for linearization}
\label{pre}

To tackle our problem's nonlinearity, we take advantage of an efficient Picard iterative scheme, as described in \cite{rpicardc, Spiridonov2019, cemnlporo,tfcmm}.  Over this section, such an iteration algorithm is presented for time-dependent multi-continuum systems.


To achieve the system \eqref{r1e}'s first goal of temporal discretization (see \cite{rpicardc, Spiridonov2019,richarde1}, for instance), we will use the following conventional backward Euler finite-difference algorithm: find $\bfa{p}_{s}=(p_{1,s},\dots, p_{N,s}) \in \bfa{V}$ such that for all $\bfa{v}=(v_1,\dots,v_N) \in \bfa{V}\,,$
\begin{align}\label{r1ed}
\left(\frac{p_{i,s+1} - p_{i,s}}{\tau} , v_i \right) + a_i(p_{i,s+1},v_i;p_{i,s+1})+ q_i(\bfa{p}_{s+1},v_i;\bfa{p}_{s+1})= (f_{i,s+1},v_i)\,, 
\end{align}
where we divide the temporal domain $[0,T]$ equally into $S$ intervals, having $\tau = T/S > 0$ as the size of time step, and 
the subscript $s$ signifies the value of a function at 
the time point $t_s= s\tau$ (with $s=0,1,\cdots,S$).

Following that, the nonlinearity in space will be linearized using Picard iteration (see \cite{rpicardc, Spiridonov2019,cemnlporo}, for instance).  At the $(s+1)$th temporal step, $\bfa{p}^0_{s+1} \in \bfa{V}$ is guessed.  Given $\bfa{p}^n_{s+1} \in \bfa{V}\,,$ with $n=0,1,2, \cdots\,,$ we seek $\bfa{p}^{n+1}_{s+1} \in \bfa{V}$ such that for all $\bfa{v}= (v_1,\dots, v_N) \in \bfa{V}\,,$ 
\begin{align}\label{r1el}
\left(\frac{p^{n+1}_{i,s+1} - p_{i,s}}{\tau} , v_i \right) + a_i(p^{n+1}_{i,s+1},v_i;p^n_{i,s+1}) + q_i(\bfa{p}^{n+1}_{s+1},v_i;\bfa{p}^n_{s+1})
= (f_{i,s+1},v_i)\,. 
\end{align}
As proved in \cite{rh2}, there exists a unique solution $\bfa{p}_{i,s+1}^{n+1}$ to this linearized system \eqref{r1el}.

When $n$ tends to $\infty\,,$ the Picard iterative procedure converges to a limit (see Appendix \ref{cp} for a theoretical proof).   
%
In simulation, we end this procedure at an $\alpha$th iteration when it satisfies a specific halting indicator, leading to the previous time data
\begin{equation}\label{pdata}
\bfa{p}_{s+1} = \bfa{p}_{s+1}^{\alpha}\,,
\end{equation}
in order to move on to the next time step in \eqref{r1ed}.  A terminating criterion is proposed over this paper employing the relative successive difference, that is, provided a user-defined tolerance $\delta_0 > 0$, if 
\begin{equation}\label{pt}
\dfrac{\|p_{i,s+1}^{n+1} - p_{i,s+1}^{n} \|_{L^2(\Omega)}}{\| p_{i,s+1}^{n} \|_{L^2(\Omega)}} \leq \delta_0\,,
\end{equation}
for $i=1,\dots,N\,,$ then the iterative process is stopped.  

Now, the fine-grid notation is considered.  First, to begin discretizing the variational problem \eqref{r1e}, we let $\mathcal{T}_h$ be a fine grid of size $h\,,$ which is assumed to be very small.  With this assumption of $h\,,$ the fine-grid solution 
will be sufficiently close to the exact solution.  Second, with respect to the rectangular fine grid $\mathcal{T}_h\,,$ we define $V_h$ as the $H^1_0(\Omega)$-conforming finite element basis space 
of piecewise bilinear functions: 
\begin{equation}\label{Vh}
V_h:= \{ v \in V: 
v |_K \in \mathcal{Q}_1(K) \; \forall K \in \mathcal{T}_h\}\,,
\end{equation}
where the space $\mathcal{Q}_1(K)$ consists of all bilinear elements (or multilinear $d$-elements when $d > 2$) over $K\,.$ We let $\bfa{V}_h = V_h^N$
and denote the $[L^2(\Omega)]^N$ projection operator onto $\bfa{V}_h$ by $P_h\,.$

On the fine scale, the completely discrete Picard iterative algorithm is as follows: beginning with an initial $\bfa{p}_{h,0} = P_h \bfa{p}_0 \in \bfa{V}_h$ having $\bfa{p}_0$ from \eqref{ini}, 
at the temporal step $(s+1)$th, we make a guess $\bfa{p}^0_{h,s+1} \in \bfa{V}_h$ 
and perform iteration from \eqref{r1el} in $\bfa{V}_h$:
\begin{align}\label{r1elh}
\begin{split}
\left(\frac{p^{n+1}_{i,h,s+1} - p_{i,h,s}}{\tau} , v_i \right) + a_i(p^{n+1}_{i,h,s+1},v_i;p^n_{i,h,s+1}) + q_i(\bfa{p}^{n+1}_{h,s+1},v_i;\bfa{p}^n_{h,s+1}) 
= (f_{i,s+1},v_i)\,,
\end{split}
\end{align}
where $n=0,1,2, \cdots\,,$ until reaching \eqref{pt} at an $\alpha$th Picard step, for all $\bfa{v}= (v_1,\dots, v_N) \in \bfa{V}_h\,.$  
To proceed to the next time step in \eqref{r1ed}, we utilize \eqref{pdata} to set the previous time data 
\begin{equation}\label{reph}
\bfa{p}_{h,s+1} = \bfa{p}_{h,s+1}^{\alpha}\,.
\end{equation}

\section{CEM-GMsFEM for coupled multi-continuum nonlinear Richards equations}\label{sec:method}

Following \cite{G2n,cemnlporo,cheung2018constraint,rtt21}, we now establish a new strategy for coupled multi-continuum nonlinear Richards equations \eqref{eq:original} in complex heterogeneous fractured porous media, using the constraint energy minimizing generalized multiscale finite element method (CEM-GMsFEM).  More specifically, in the pressure computation for the equivalently nonlinear system \eqref{r1e}, we will show the establishment of auxiliary space (utilizing a novel idea of two sample sources) and multiscale space.  To properly construct such CEM-GMsFEM, after temporal discretization of \eqref{r1e}, the linearized formulation \eqref{r1el} can be used to consider the nonlinearity as a constant at each Picard iterative step. Multiscale space can therefore be built according to this nonlinearity.

\subsection{Overview}\label{over}

\begin{figure}[ht!]
\centering
\includegraphics[width=0.5\linewidth]{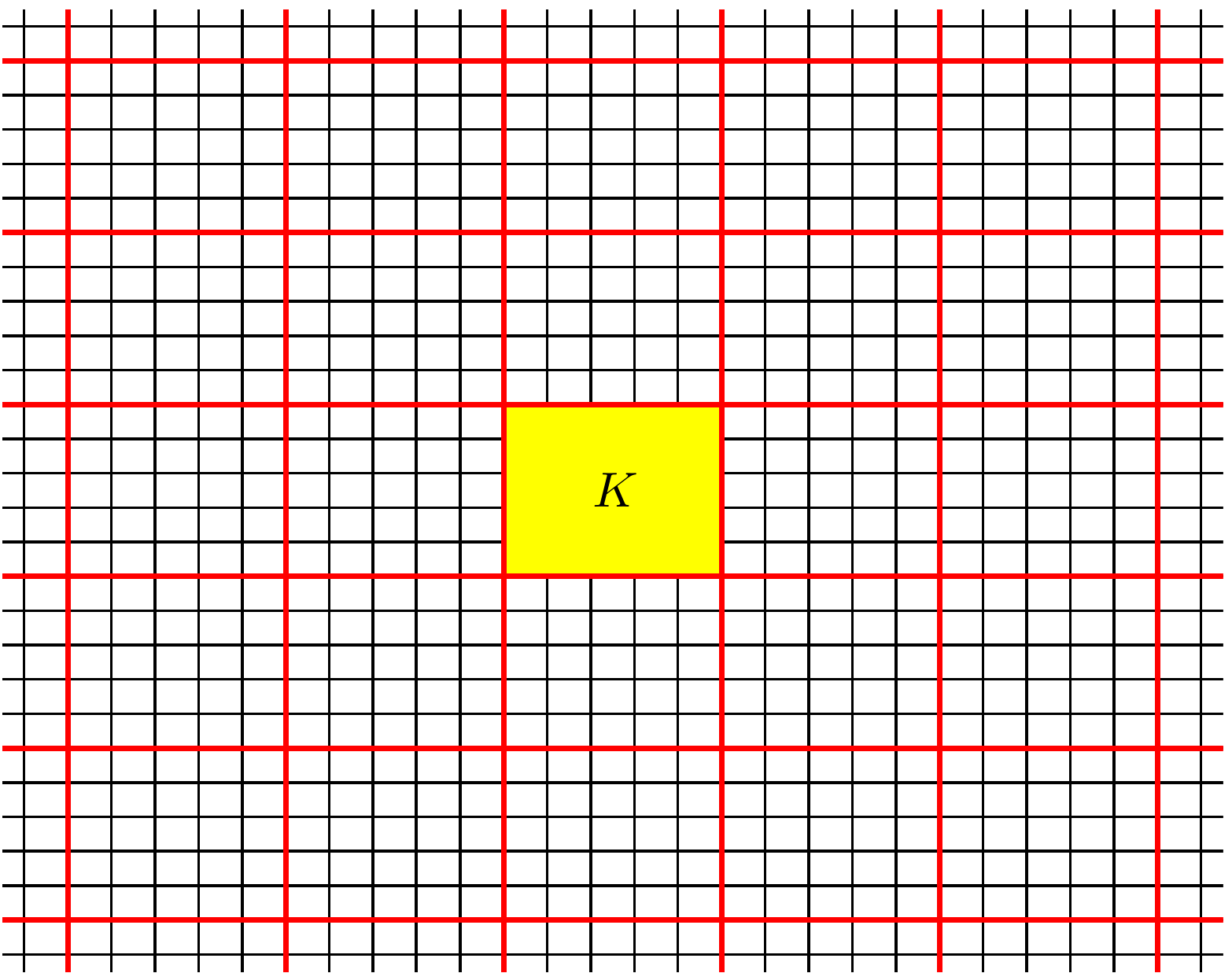}
\caption{Illustration of the fine grid, coarse grid, and a coarse block $K\,.$}
\label{fig:mesh}
\end{figure}

First, we will go over the concepts of coarse and fine grids.  The start is partitioning $\Omega$ into finite elements, where multiscale characteristics are not necessarily resolved. 
The partition is named coarse grid $\mathcal{T}^H\,,$ and $\mathcal{T}_h$ is its refinement.  In $\mathcal{T}^H\,,$ a generic element $K$ is named a coarse-grid block (also known as coarse element or coarse patch). 
Moreover, we call $H > 0$ the coarse-grid size, where $H \gg h\,.$
Let $N_c$ be the number of coarse blocks and $N_v$ be the number of coarse-grid nodes.  The collection of all coarse nodes (vertices) is denoted by $\{\bfa{x}_k\}_{k=1}^{N_v}\,.$  
Figure~\ref{fig:mesh} illustrates the fine grid, coarse grid, as well as a coarse block $K\,.$  

For some subdomain $D \subset \Omega\,,$ the restrictions of $V$ and $\bfa{V}$ on $D$ are respectively denoted by $V(D)$ and $\bfa{V}(D)\,.$  Furthermore, the subspace of $\bfa{V}(D)$ containing functions with zero trace on $\partial D$ is represented by $\bfa{V}_0(D)\,.$  Using this definition on a coarse block $K_j\,,$ provided $\bfa{u} = (u_1, \cdots, u_N) \in \bfa{V}(K_j)\,,$ for all $\bfa{p} = (p_1, \cdots, p_N)\,, \ \bfa{v} = (v_1, \cdots, v_N) \in \bfa{V}(K_j)\,, i = 1,2, \cdots,N\,,$ we define local bilinear forms by
\begin{equation}
\begin{aligned}
a^{(j)}_i(p_i, v_i; u_i) & = 
\int_{K_j} \varkappa_i (u_i) \nabla p_i \cdot \nabla v_i \, dx\,, 
\\
q^{(j)}_i(\bfa{p},v_i; \bfa{u}) & = \sum_{l=1}^N \int_{K_j} Q_{il}(u_i,u_l)(p_i-p_{l}) v_i \, \dx\,, 
\\
r^{(j)}_i(p_i, v_i; u_i) & = 
\int_{K_j} \widetilde{\varkappa}_i (u_i) p_i  v_i \, dx 
\,.
\end{aligned}
\label{eq:bilinear_loc1}
\end{equation}
Here, 
\[\widetilde{\varkappa}_i(u_i) = \varkappa_i(u_i) \sum_{k=1}^{N_v} \vert \nabla \chi_{k,i} \vert^2\,,\] 
in which for the coarse node $\bfa{x}_k\,,$ each $\chi_{k,i}$ (with linear boundary
conditions for cell problems \cite{pou}) is a conventional multiscale finite element basis function within the $i$th continuum.  In the current context, $\{\chi_{k,i}\}_{k=1}^{N_v}$ is a collection of bilinear partition of unity functions (for $\mathcal{T}^H$) supported in the continuum $i$th.  On each coarse block $K_j\,,$ the coupled local bilinear forms are defined as follows:
for any $\bfa{p} = (p_1, \cdots, p_N), \bfa{v} = (v_1, \cdots, v_N) \in \bfa{V}(K_j)$, $i = 1, \cdots,N\,,$ 
\begin{equation}
\begin{aligned}
a^{(j)}(\bfa{p},\bfa{v}; \bfa{u}) & = \sum_i a^{(j)}_i(p_i, v_i; u_i), \\
q^{(j)}(\bfa{p},\bfa{v}; \bfa{u}) & = \sum_i q^{(j)}_i(\bfa{p}, v_i; \bfa{u}), \\
r^{(j)}(\bfa{p},\bfa{v}; \bfa{u}) & = \sum_i r^{(j)}_i(p_i, v_i; u_i), \\
a^{(j)}_Q(\bfa{p},\bfa{v}; \bfa{u}) & = a^{(j)}(\bfa{p},\bfa{v}; \bfa{u}) + q^{(j)}(\bfa{p},\bfa{v}; \bfa{u}).
\end{aligned}
\label{eq:bilinear_loc}
\end{equation}
We also define the bilinear forms $a_Q$ and $r$ by
\begin{align}\label{bis}
\begin{split}
a_Q(\bfa{p},\bfa{v}; \bfa{u}) &= \sum_{j=1}^{N_c} a^{(j)}_Q(\bfa{p},\bfa{v}; \bfa{u})\,, \\
r(\bfa{p},\bfa{v};\bfa{u}) &= \sum_{j=1}^{N_c} r^{(j)}(\bfa{p},\bfa{v};\bfa{u})\,.
\end{split}
\end{align}

In multiscale space, our major target is seeking for \eqref{r1e} a multiscale solution $\bfa{p}_{\tu{ms}}$ that approximates the fine-scale solution $\bfa{p}_h$ better than using the GMsFEM (\cite{gnone}).   For this purpose, the CEM-GMsFEM is utilized to attain the multiscale solution $\bfa{p}_{\tu{ms}}\,.$  Two levels are required to build the multiscale space.  First, through the GMsFEM, an auxiliary space is generated.   Second, a multiscale space $\bfa{V}_{\tu{ms}}$ is established (benefiting from that auxiliary space) and possesses multiscale basis functions with locally minimal energy over some subregions.
Ultimately, a multiscale solution can be found using these 
multiscale basis functions.  Note that the obtained $\bfa{V}_{\tu{ms}}$ is stable throughout this procedure when utilizing either the first or second source of samples $\bfa{u}_s$ from Subsection \ref{auxbf}.  We refer the readers to \cite{rtt21, mcl, G1, G2n, chungres1, chung2015residual, chung2016adaptive,gne}
and \cite{cem1,cem2f,cheung2018constraint, cemnlporo} for more information about the GMsFEM and CEM-GMsFEM, respectively.  At the first stage, the GMsFEM will be used to design our auxiliary multiscale basis functions as follows. 

\subsection{Auxiliary multiscale basis functions}\label{auxbf}
Given a set of samples $\{\bfa{u}_b\}_{b=1}^B 
\subset \bfa{V}$ 
and weights $\{w_b\}_{b=1}^B \subset \mathbb{R}^+\,.$  The following bilinear forms are determined from weighted Monte Carlo integration \cite{samples}: 
\begin{equation}\label{lbi}
	\begin{split}
		A_Q^{(j)}(\bfa{p},\bfa{v}) & = \sum_{b=1}^{B} w_b\, a_Q^{(j)}(\bfa{p},\bfa{v}; \bfa{u}_b), \\
		R^{(j)}(\bfa{p},\bfa{v}) & = \sum_{b=1}^{B} w_b \, r^{(j)}(\bfa{p},\bfa{v}; \bfa{u}_b).
	\end{split}
\end{equation}


\bigskip

Next, we will create our local auxiliary multiscale basis functions employing the GMsFEM.  In particular, these coupled functions are identified by a local spectral problem, that is, to seek a real number $\lambda_k^{(j)}$ and an associated function $\bfa{\phi}_k^{(j)} \in \bfa{V}(K_j)$ such that
\begin{equation}\label{eq:spectral_prob}
A_Q^{(j)}(\bfa{\phi}_k^{(j)},\bfa{v}) = \lambda_k^{(j)} R^{(j)}(\bfa{\phi}_k^{(j)},\bfa{v}) \text{ for any } \bfa{v} \in \bfa{V}(K_j)\,.
\end{equation}
{The eigenfunctions $\bfa{\phi}_k^{(j)}$ of \eqref{eq:spectral_prob} are normalized in the norm produced by the inner product $R^{(j)}$ as follows:
\begin{equation}\label{normalize}
R^{(j)}(\bfa{\phi}_k^{(j)},\bfa{\phi}_k^{(j)}) = 1\,.
\end{equation}
The eigenvalues $\lambda_k^{(j)}$ of \eqref{eq:spectral_prob} are organized in nondecreasing order over $k\,.$  Then, using the first corresponding $L_j$ eigenfunctions, we generate the following local auxiliary multiscale space:
\begin{equation}
\bfa{V}_{\tu{aux}}^{(j)} =
\bfa{V}_{\tu{aux}}(K_j)
= \text{span} \{ \bfa{\phi}_k^{(j)}: 1 \leq k \leq L_j\}\,.
\end{equation}
The sum of such local auxiliary multiscale spaces represents the global auxiliary multiscale space: 
\begin{equation}
\bfa{V}_{\tu{aux}} = \oplus_{j=1}^{N_c} \bfa{V}_{\tu{aux}}^{(j)}\,.
\end{equation}
Also, the global bilinear forms $A_Q$ and $R$ are defined by
\begin{align}\label{biAS}
	\begin{split}
		A_Q(\bfa{p},\bfa{v}) &= \sum_{j=1}^{N_c} A^{(j)}_Q(\bfa{p},\bfa{v})\,, \\
		R(\bfa{p},\bfa{v}) &= \sum_{j=1}^{N_c} R^{(j)}(\bfa{p},\bfa{v})\,.\\
    \end{split}
\end{align}

\bigskip

We briefly mention two sources of samples $\{\bfa{u}_b\}_{b=1}^B$ which will be considered in our numerical experiments. 
The first source of samples is simply from choosing 
a single realization $\bfa{p}_h^\star \in \bfa{V}_h$, that is, $B = 1$ where $w_1 = 1$ in \eqref{lbi}, and $\bfa{u}_1 = \bfa{p}_h^\star$ is the steady-state FEM solution \eqref{reph}.
This first sample source will be used for steady-state cases.

The second source of samples is from taking realizations of the sink or source function in \eqref{eq:original} and 
the backward Euler temporal discretization \eqref{r1ed} within the fully Picard discrete scheme \eqref{r1elh} for solving \eqref{eq:original}, 
to obtain from \eqref{reph} a numerical fine-grid approximation $\bfa{p}^\star_{h,s} \in \bfa{V}_h$ at the time step $s$ (where $t = s\tau$).
In this case, we have $B = S+1$, and 
\begin{equation}\label{seconds}
	\bfa{u}_b =  \bfa{p}^\star_{h,{b-1}}\,.
\end{equation}  
By setting the weights 
\begin{equation}
	w_b = \begin{cases}
		1 & \text{ if } b = 2,3,\ldots,S\,, \\
		1/2 & \text{ if } b = 1,S+1\,,
	\end{cases}
\end{equation}
the weighted integration \eqref{lbi} can be regarded as 
numerical integration of piecewise linear functions using trapezoidal rule in the temporal variable.  That is, 
\begin{equation}\label{lag}
	\begin{split}
		A_Q^{(j)}(\bfa{p},\bfa{v}) & = \int_0^T  a_Q^{(j)}(\bfa{p},\bfa{v}; \bfa{p}_h^\star)\, \dt \,, \\
		R^{(j)}(\bfa{p},\bfa{v}) & = \int_0^T  r^{(j)}(\bfa{p},\bfa{v}; \bfa{p}_h^\star) \, \dt, 
	\end{split}
\end{equation}
where $\bfa{p}_h^\star$ is the piecewise linear Lagrange interpolation of $\{\bfa{p}_{h,s}^\star\}_{s=0}^S$ on the temporal grid $\{s\tau\}_{s=0}^S$. 
This second sample source will be applied to the time-dependent cases. In general, one can utilize multiple realizations in a similar manner for both sample sources. 


\subsection{Multiscale space}\label{msp}

To define our multiscale basis functions for spanning the solution space, we introduce the concept of $\phi$-orthogonality. Providing an auxiliary basis function $\bfa{\phi}_k^{(j)} \in \bfa{V}_{\tu{aux}}$ within a coarse block $K_j\,,$ we state that $\bfa{\psi} \in \bfa{V}$ is $\bfa{\phi}_k^{(j)}$-orthogonal if for $1 \leq k' \leq L_{j'}$ and $1 \leq j' \leq N_c\,,$ 
\begin{equation}
R \left(\bfa{\psi}, \bfa{\phi}_{k'}^{(j')}\right) =  \delta_{k,k'} \delta_{j,j'}\,,
\end{equation}
equivalently,
\begin{equation}
R \left(\bfa{\psi}, \bfa{\phi}_{k}^{(j)}\right) = 1\,, \quad 
R \left(\bfa{\psi}, \bfa{\phi}_{k'}^{(j')}\right) = 0
\quad \tu{ with } k' \neq k \tu{ or } j' \neq j \,. 
\end{equation}
This $\bfa{\phi}_k^{(j)}$-orthogonality gives rise to the orthogonal projection operator $\pi: [L^2(\Omega)]^N \to \bfa{V}_{\tu{aux}}$ proposed by 
$\dd \pi = \sum_{j=1}^{N_c} \pi_j$, where $\pi_j: [L^2(K_j)]^N \to \bfa{V}_{\tu{aux}}^{(j)}$ is defined as
\begin{equation}
\pi_j(\bfa{v}) = \sum_{k=1}^{L_j} \dfrac{R^{(j)}(\bfa{v},\bfa{\phi}_k^{(j)})}{R^{(j)}(\bfa{\phi}_k^{(j)},\bfa{\phi}_k^{(j)})} \bfa{\phi}_k^{(j)} \text{ for any } \bfa{v} \in [L^2(K_j)]^N\,.
\end{equation}

Our global multiscale basis functions are now being built.  For each auxiliary function $\bfa{\phi}_k^{(j)} \in \bfa{V}_{\tu{aux}}\,,$ the solution to the following constrained energy minimization problem specifies the global multiscale basis function $\bfa{\psi}_{k}^{(j)} \in \bfa{V}\,:$ 
\begin{equation}
\bfa{\psi}_{k}^{(j)} = \argmin \left\{ A_Q(\bfa{\psi}, \bfa{\psi}) : 
\bfa{\psi} \in \bfa{V} \text{ is } \bfa{\phi}_k^{(j)} \text{-orthogonal}\right\}.
\label{eq:min1_glo}
\end{equation}
The variational form of this minimization problem \eqref{eq:min1_glo} is as follows: 
determine $\bfa{\psi}_{k}^{(j)} \in \bfa{V}$ and $\bfa{\mu}_{k}^{(j)} \in \bfa{V}_{\tu{aux}}$ such that
\begin{equation}
\begin{split}
A_Q(\bfa{\psi}_{k}^{(j)}, \bfa{v}) + R(\bfa{v}, \bfa{\mu}_{k}^{(j)}) & = 0 \text{ for all } \bfa{v} \in \bfa{V}, \\
R(\bfa{\psi}_{k}^{(j)} - \bfa{\phi}_k^{(j)}, \bfa{\nu}) & = 0 \text{ for all } \bfa{\nu} \in \bfa{V}_{\tu{aux}}.
\end{split}
\label{eq:var1_glo}
\end{equation}

We identify our localized multiscale basis functions as a result of the establishment of global multiscale basis functions.
An oversampled domain is created by extending the coarse grid block $K_{{j}}$ by $m$ coarse-grid layers, for each $K_{{j}} \in \mathcal{T}^H\,.$  Figure~\ref{fig:oversample} depicts an example of an oversampled region.  The solution of 
the following constrained energy minimization problem defines the localized multiscale basis function $\bfa{\psi}_{k,{\tu{ms}}}^{(j)} \in \bfa{V}_0(K_{j,m})\,:$
\begin{equation}
\bfa{\psi}_{k,{\tu{ms}}}^{(j)} = \argmin \left\{ A_Q(\bfa{\psi}, \bfa{\psi}) : 
\bfa{\psi} \in V_0(K_{j,m}) \text{ is } \bfa{\phi}_k^{(j)} \text{-orthogonal}\right\}.
\label{eq:min1}
\end{equation}
Now, let $Z_{j}$ be the set of indices such that if $z \in Z_j\,,$ then $K_{z} \subset K_{j,m}\,;$ 
and let
\[\bfa{V}_{\tu{aux}}(K_{j,m}) = \sum_{z \in Z_j} \bfa{V}_{\tu{aux}}(K_{z}) = \tu{span} \left\{ \bfa{\phi}_k^{(z)} : 1 \leq k \leq L_{z}\,, z \in Z_j \right \}
\,.\]
Then, the following variational problem is equivalent to the minimization problem \eqref{eq:min1}: 
find $\bfa{\psi}_{k,{\tu{ms}}}^{(j)} \in \bfa{V}_0(K_{j,m})$ and $\bfa{\mu}_{k,\tu{ms}}^{(j)} \in \bfa{V}_{\tu{aux}}(K_{j,m})$ 
such that
\begin{equation}
\begin{split}
A_Q(\bfa{\psi}_{k,{\tu{ms}}}^{(j)}, \bfa{v}) + R(\bfa{v}, \bfa{\mu}_{k,\tu{ms}}^{(j)}) & = 0 \text{ for all } \bfa{v} \in \bfa{V}_0(K_{j,m}), \\
R(\bfa{\psi}_{k,{\tu{ms}}}^{(j)} - \bfa{\phi}_k^{(j)}, \bfa{\nu}) & = 0 \text{ for all } \bfa{\nu} \in \bfa{V}_{\tu{aux}}(K_{j,m})\,.
\end{split}
\label{eq:var1}
\end{equation}

\begin{figure}[ht!]
\centering
\includegraphics[width=0.5\linewidth]{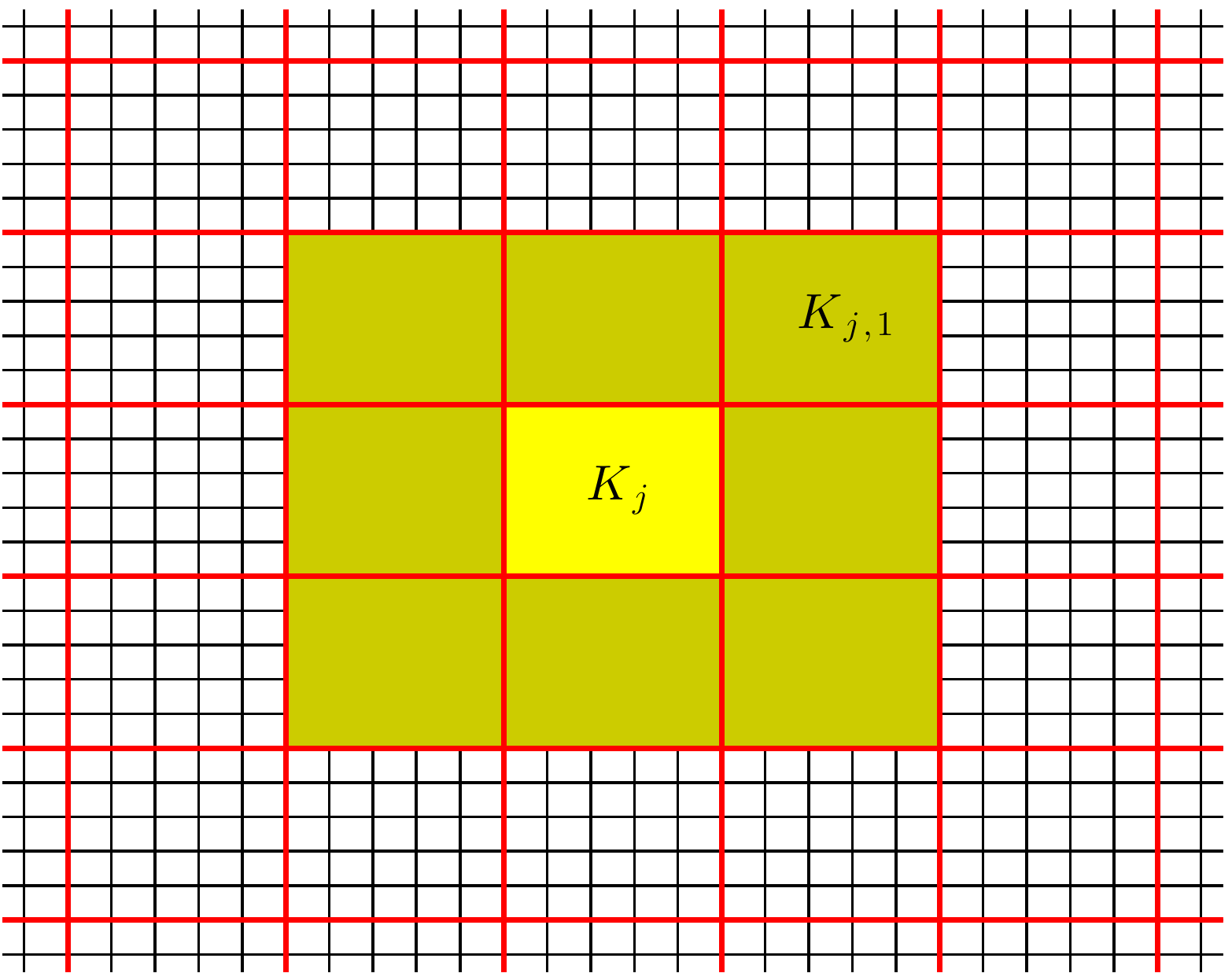}
\caption{Illustration of an oversampled region $K_{j,1}$ established via extending $K_j$ by $1$ coarse-grid layer.}
\label{fig:oversample}
\end{figure}

Using the localized multiscale basis functions, we generate
the multiscale finite element space
\begin{equation}\label{vms}
\bfa{V}_{\tu{ms}} = \text{span}  \{\bfa{\psi}_{k,{\tu{ms}}}^{(j)} : 1 \leq k \leq L_j, 1 \leq j \leq N_c \}\,.
\end{equation}
Remark that outside of some local (oversampled) subdomains, the global multiscale basis functions \eqref{eq:min1_glo} (which are generally supported in the whole domain) exponentially decay \cite{cem1}.
This characteristic is critical in the CEM-GMsFEM's convergence investigation, demonstrating the utilization of local multiscale basis functions \eqref{eq:min1} in $\bfa{V}_{\tu{ms}}$ \cite{cheung2018constraint}.

\subsection{CEM-GMsFEM for coupled system of nonlinear Richard equations}\label{alg}

We note first that the space $\bfa{V}=\bfa{H}_0^1(\Omega)$ is continuous throughout the previous Sections and Subsections \ref{over}--\ref{msp}. 
The multiscale space $\bfa{V}_{\tu{\tu{ms}}} \subset \bfa{V}$ therefore requires some finite dimensional analogues for simulations.  As a result, the problem under consideration is tackled via the fine grid together with a suitable finite element method \cite{cem2f,cemnlporo} in our numerical computations.

Recalling $\bfa{p}_{h,0} = P_h \bfa{p}_0 \in \bfa{V}_h$ with $\bfa{p}_0$ from \eqref{ini}, we now have an initial $\bfa{p}_{\tu{ms}}(0, \cdot)= \bfa{p}_{\tu{ms},0}$ that satisfies 
\begin{equation}\label{pms0}
	a_i(p_{i,h,0} - p_{i,\tu{ms},0},v_i;p_{i,h,0}) + q_i(\bfa{p}_{h,0}-\bfa{p}_{\tu{ms},0},v_i;\bfa{p}_{h,0}) = 0\,,
\end{equation}
for any $\bfa{v} = (v_1, \cdots,v_N) \in \bfa{V}_{\tu{ms}}\,.$ 
Fixing the $(s+1)$th temporal step, our strategy (as in \cite{rtt21, gne, cemnlporo}) is to solve the problem 
(\ref{r1e}) through linearization relied on Picard's iterative technique.  This can be done by employing at each Picard iteration the CEM-GMsFEM (in Subsections \ref{auxbf} and \ref{msp}) with the constructed stable offline multiscale space $\bfa{V}_{\tu{ms}}$ (proposed at the close of Subsection \ref{over}). 

In particular, all along the online stage, the full model reduction approach is as follows: beginning with $\bfa{p}_{\tu{ms},0} \in \bfa{V}_{\tu{ms}}$ from \eqref{pms0}, we pick a guess $\bfa{p}^0_{\tu{ms},s+1}\in \bfa{V}_{\tu{ms}}$ at the temporal step $(s+1)$th 
and iterate from \eqref{r1el} in $\bfa{V}_{\tu{ms}}$:
\begin{align}\label{ongmspicard}
\left(\frac{p^{n+1}_{i,\tu{ms},s+1} - p_{i,\tu{ms},s}}{\tau} , \psi_i \right) + a_i(p^{n+1}_{i,\tu{ms},s+1},\psi_i;p^n_{i,\tu{ms},s+1}) + q_i(\bfa{p}^{n+1}_{\tu{ms},s+1},\psi_i;\bfa{p}^n_{\tu{ms},s+1}) = (f_{i,s+1},\psi_i)\,,
\end{align}
where $\bfa{\psi} = (\psi_1,\cdots, \psi_N)\in \bfa{V}_{\tu{ms}}$ and $n=0,1,2, \cdots\,,$ until getting to \eqref{pt} at some Picard step $\alpha$th.  
To move to the next time step in \eqref{r1ed}, we use \eqref{pdata} for selecting the previous temporal data 
\begin{equation}\label{predatacem}
\bfa{p}_{\tu{ms},s+1} = \bfa{p}_{\tu{ms},1}^{\alpha}\,.
\end{equation}


\bigskip

\begin{remark}\label{cemgen}
The CEM-GMsFEM for single-continuum cases are treated similarly to the multi-continuum cases in this Section by allowing the transfer term to vanish in \eqref{r1e} and in all related expressions.  Moreover, the steady-state cases are handled similarly to the time-dependent cases in this Section by letting the second source of samples \eqref{seconds} for Subsection \ref{auxbf} to be the first source of samples  \eqref{reph} 
as the steady-state FEM solution.  After we choose an appropriate source of samples,  the snapshot functions as well as the basis functions are time-independent.  Relying on specific case, Eq. \eqref{ongmspicard} is the corresponding equation in Section \ref{pre} (time-dependent dual-continuum case \eqref{r1elh}, time-dependent single-continuum case, 
steady-state dual-continum case,
and steady-state single-continum case),
where the subscript ``$h$'' is replaced by ``$\tu{ms}$''.

\end{remark}

\section{Numerical examples}\label{sec:numerical}

We will give various numerical tests in this section to demonstrate our approach's performance.  In each test, the effect of coarse-grid size $H$ is investigated. For all experiments, we used the fine-grid size $h = 1/128$ and the number of oversampling layers $m \approx 10 \log(1/H) / \log(64)$. Within each experiment, the number of local multiscale basis functions is fixed throughout all coarse elements. 
We compare the solutions obtained by our strategy employing the constraint energy minimizing
generalized multiscale finite element method (CEM-GMsFEM, abbreviated by CEM)  
with the solutions computed by the finite element method (FEM). 

The spatial domain is $\Omega = [0,1]^2\,.$  For time-dependent equations in the temporal interval $[0,T]\,,$ the stopping time is $T=S\tau=2$ while the temporal step size is $\tau = 1/10\,,$ so there are $S=20$ time steps.  We consider the channelized media and consequently deal with the high-contrast coefficients defined in the spatial domain $\Omega\,.$ 
The Picard iteration's halting indicator is $\delta_0 = 10^{-5}$, which guarantees the convergence of this linearization procedure.
The continua are assumed to be isotropic (and the anisotropic case is treated in the same way).  Then, hydraulic conductivity tensors can be considered as scalar functions $\varkappa_i(\bfa{x},p_i)$ (multiplying with the identity matrix \cite{kanhe}), for $i=1,2\,.$ 


\subsection{Experiements for single-continuum Richards equations}\label{1con}

In this section, we examine steady-state and time-dependent single-continuum Richards equations from \eqref{eq:original}, respectively: find $p \in V$ such that
\beq
\label{eq:original0}
- \div \left(\kappa(\bfa{x}) e^{p(\bfa{x})} \nabla p(\bfa{x}) \right) = 1 \  \textrm{in } \Omega\,,
\eeq
and
\begin{equation}
\label{eq:original00}
\frac{\partial p(t,\bfa{x})}{\partial t} - \div [\kappa(\bfa{x}) e^{p(t,\bfa{x})}\nabla p(t,\bfa{x})]
= f(t,\bfa{x}) \  \textrm{in} \ (0,T] \times \Omega\,,
\end{equation}
where $f(t,\bfa{x}) = \sin(\pi x_1)\sin(\pi x_2)\,,$ and $\varkappa(\bfa{x},p) = \kappa(\bfa{x}) e^{p(t,\bfa{x})}\,.$  Both problems have the zero Dirichlet boundary condition, and the initial condition for (\ref{eq:original00}) is $p(0,\bfa{x})=0\,.$  The permeability field $\kappa(\bfa{x})$ in (\ref{eq:original0}) and (\ref{eq:original00}) is depicted in Figure \ref{fig:func_a_single}. The value of $\kappa(\bfa{x})$ in the yellow regions (channels) is $1000$ and in the blue region is $10\,.$
\begin{figure}[ht!]
	\centering
	\includegraphics[width=0.6\linewidth]{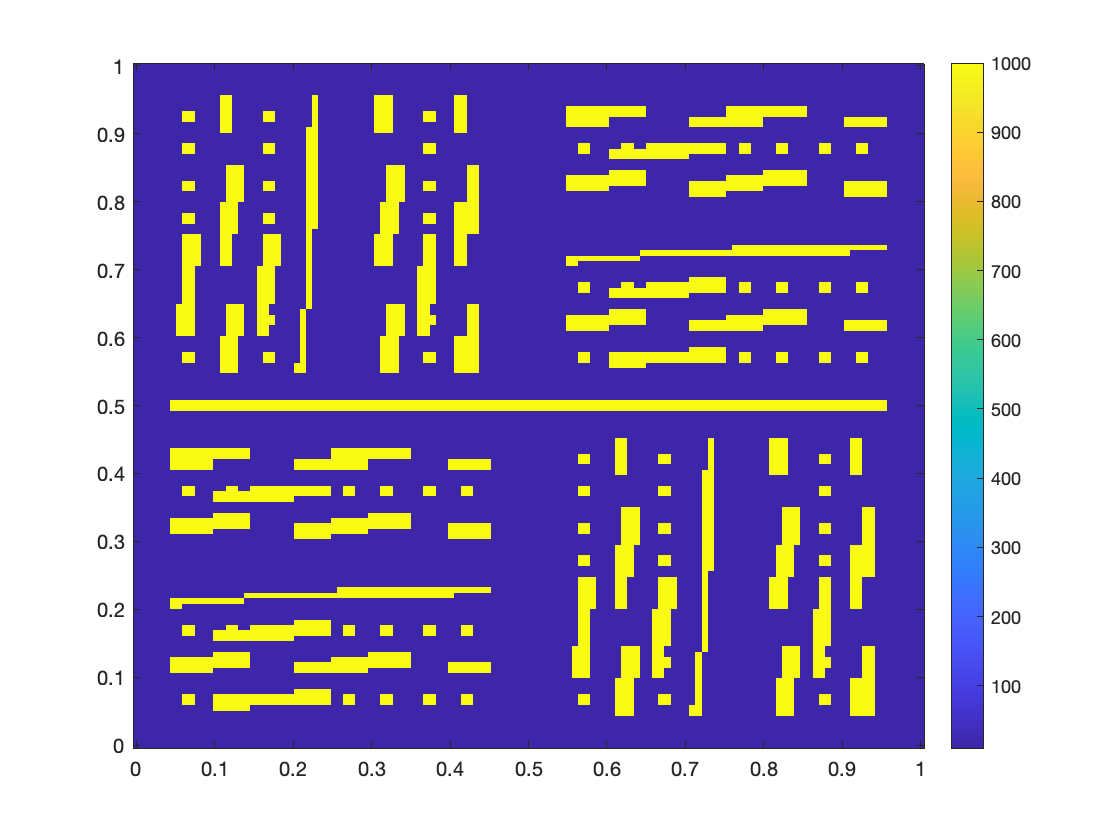}
	\caption{Permeability field $\kappa(\bfa{x})$ for problems (\ref{eq:original0}) and (\ref{eq:original00}).}
	\label{fig:func_a_single}
\end{figure}

We compute the CEM solutions $p_{\tu{ms}}$ from \eqref{ongmspicard} as well as their respective FEM references $p_h$ in \eqref{r1elh} for Eqs.\ (\ref{eq:original0}) and 
(\ref{eq:original00}). 
Note that the above single-continuum equations do not have the interaction terms, thus the corresponding terms in (\ref{r1elh}) and \eqref{ongmspicard} are ignored when following the algorithms to compute the CEM and FEM solutions. Further, for steady-state case (\ref{eq:original0}), we do not need to involve the iteration with respect to time.  The initial guess for Picard iteration of problem \eqref{eq:original0} is the identically zero function in the domain and of problem (\ref{eq:original00}) is the previous time data (\eqref{reph} for FEM and \eqref{predatacem} for CEM). To construct multiscale basis functions, the sampling method described in Subsection \ref{auxbf} is utilized. We employed the second source of samples \eqref{seconds} based on \eqref{reph} for time-dependent case and used the steady-state FEM solution $p_h$ \eqref{reph} as the first source of samples for the steady-state case. 

In both experiments, the relative $L^2$ and $H^1$ errors between our CEM solutions $p_{\tu{ms}}$ 
and the FEM references $p_{h}$ 
are defined as 
\begin{align}\label{errors_single}
	\begin{split}
e^p_{L^2}=\frac{||p_{\tu{ms}}-p_{h}||_{L^2(\Omega)}}{||p_{h}||_{L^2(\Omega)}}, \
e^p_{H^1}=\frac{||\nabla p_{\tu{ms}}-\nabla p_{h}||_{\bfa{L}^2(\Omega)}}{||\nabla p_{h}||_{\bfa{L}^2(\Omega)}}\,. 
\end{split}
\end{align}
We investigate these errors with respect to the coarse-grid size $H$, by different number of local multiscale basis functions. 
Table \ref{tab:steady_single} presents the errors for the steady-state equation (\ref{eq:original0}), and Table \ref{tab:td_single} shows the errors in the time-dependent case (\ref{eq:original00}). We note that the total number of degrees of freedom for our multiscale method ($\textrm{dim}(V_{\tu{ms}})$) relates entirely to the coarse-grid size $H$ and the number of local multiscale basis functions.
In all tables, we observe that the numerical approximations are very accurate for every choice of coarse-grid size $H$. Also, $H^1$ errors less than $3\%$ and $L^2$ errors less than $0.5\%$ even with relatively large coarse-grid size $H=1/4$, where only $96$ degrees of freedom are utilized for CEM at maximum throughout the experiments. This number is much less than $16129$, the total number of degrees of freedom used in FEM.  It is explicit from those tables that as the sequence of coarse-grid sizes $H$ converges, the sequence of CEM solutions converges.
According to the tables, both $H^1$ and $L^2$ errors can be further decreased once more local multiscale basis functions and oversampling layers are involved.  However, a too large number of multiscale basis functions has a direct impact on the method's computational complexity. It is unknown whether the contrast has a direct effect on the number of required multiscale basis functions. 
Figure \ref{fig:sol_tdsingle} illustrates the comparison of the plots of solutions to (\ref{eq:original00}) at the final time $T=2$, computed by the CEM and FEM when $H=1/16$. One can see that each solution obtained by CEM almost coincides with its reference solution computed by FEM.

\begin{table}[ht!]
\centering
\begin{subtable}{0.8\textwidth}
\centering
\begin{tabular}{|c|c|c||c||c|}
\hline
$H$ & $m$ &$\textrm{dim}(V_{\tu{ms}})$ &$H^1$ error & $L^2$ error\\
\hline
$1/4$ & 3 &64& 2.0760\% & 0.3446\% \\
$1/8$ & 5 &256& 1.3022\% & 0.1188\% \\
$1/16$ & 7 &1024&0.7354\% & 0.0389\% \\
$1/32$ & 8 &4096& 0.2829\% &  0.0073\% \\
\hline
\end{tabular}
\caption{Errors with $4$ local basis functions for (\ref{eq:original0}).}
\label{tab:steady_single_J4}
\end{subtable}
\begin{subtable}{0.8\textwidth}
\centering
\begin{tabular}{|c|c|c||c||c|}
\hline
$H$ & $m$&$\textrm{dim}(V_{\tu{ms}})$ & $H^1$ error & $L^2$ error \\
\hline
$1/4$ & 3 &80& 1.6471\% & 0.2410\% \\
$1/8$ & 5 &320& 0.9370\%  & 0.0885\% \\
$1/16$ & 7 &1280&0.4682\% & 0.0219\%  \\
$1/32$ & 8 &5120& 0.1736\%  & 0.0041\%  \\
\hline
\end{tabular}
\caption{Errors with $5$ local basis functions for (\ref{eq:original0}).}
\label{tab:td_single_J5}
\end{subtable}
\begin{subtable}{0.8\textwidth}
\centering
\begin{tabular}{|c|c|c||c||c|}
\hline
$H$ & $m$&$\textrm{dim}(V_{\tu{ms}})$ & $H^1$ error & $L^2$ error \\
\hline
$1/4$ & 3 &96& 1.4018\% & 0.1914\% \\
$1/8$ & 5 &384& 0.6085\%  & 0.0452\% \\
$1/16$ & 7 &1536&0.2779\% & 0.0112\%  \\
$1/32$ & 8 &6144& 0.0707\%  & 0.0012\%  \\
\hline
\end{tabular}
\caption{Errors with $6$ local basis functions for (\ref{eq:original0}).}
\label{tab:steady_single_J6}
\end{subtable}
\vspace*{-3mm}
\caption{Relative $L^2$ and $H^1$ errors for (\ref{eq:original0}) with different number of local basis functions; $\textrm{dim}(V_h) = 16129\,.$}
\label{tab:steady_single}
\end{table}


\begin{table}[ht!]
\centering
\begin{subtable}{0.8\textwidth}
\centering
\begin{tabular}{|c|c|c||c||c|}
\hline
$H$ & $m$&$\textrm{dim}(V_{\tu{ms}})$ & $H^1$ error & $L^2$ error\\
\hline
$1/4$ & 3 &64& 2.7544\% & 0.4359\% \\
$1/8$ & 5 &256& 1.3024\% & 0.1261\% \\
$1/16$ & 7 &1024&0.7187\% & 0.0376\% \\
$1/32$ & 8 &4096& 0.2687\% &  0.0068\% \\
\hline
\end{tabular}
\caption{Errors with $4$ local basis functions for (\ref{eq:original00}).}
\label{tab:td_single_J4}
\end{subtable}
\begin{subtable}{0.8\textwidth}
\centering
\begin{tabular}{|c|c|c||c||c|}
\hline
$H$ & $m$ &$\textrm{dim}(V_{\tu{ms}})$& $H^1$ error & $L^2$ error \\
\hline
$1/4$ & 3 &80& 2.0580\% & 0.2887\% \\
$1/8$ & 5 &320& 1.0094\%  & 0.0962\% \\
$1/16$ & 7 &1280&0.4816\% & 0.0226\%  \\
$1/32$ & 8 &5120& 0.1659\%  & 0.0038\%  \\
\hline
\end{tabular}
\caption{Errors with $5$ local basis functions for (\ref{eq:original00}).}
\label{tab2:td_single_J5}
\end{subtable}
\begin{subtable}{0.8\textwidth}
\centering
\begin{tabular}{|c|c|c||c||c|}
\hline
$H$ & $m$ &$\textrm{dim}(V_{\tu{ms}})$& $H^1$ error & $L^2$ error \\
\hline
$1/4$ & 3 &96& 1.8294\% & 0.2391\% \\
$1/8$ & 5 &384& 0.7953\%  & 0.0676\% \\
$1/16$ & 7 &1536&0.2685\% & 0.0106\%  \\
$1/32$ & 8 &6144& 0.0701\%  & 0.0011\%  \\
\hline
\end{tabular}
\caption{Errors with $6$ local basis functions for (\ref{eq:original00}).}
\label{tab:td_single_J6}
\end{subtable}
\vspace*{-3mm}
\caption{Relative $L^2$ and $H^1$ errors for (\ref{eq:original00}) with different number of local basis functions; $\textrm{dim}(V_h) = 16129\,.$}
\label{tab:td_single}
\end{table}


\begin{figure}[ht!]
\centering
\begin{subfigure}{0.45\textwidth}
  \includegraphics[width=\textwidth]{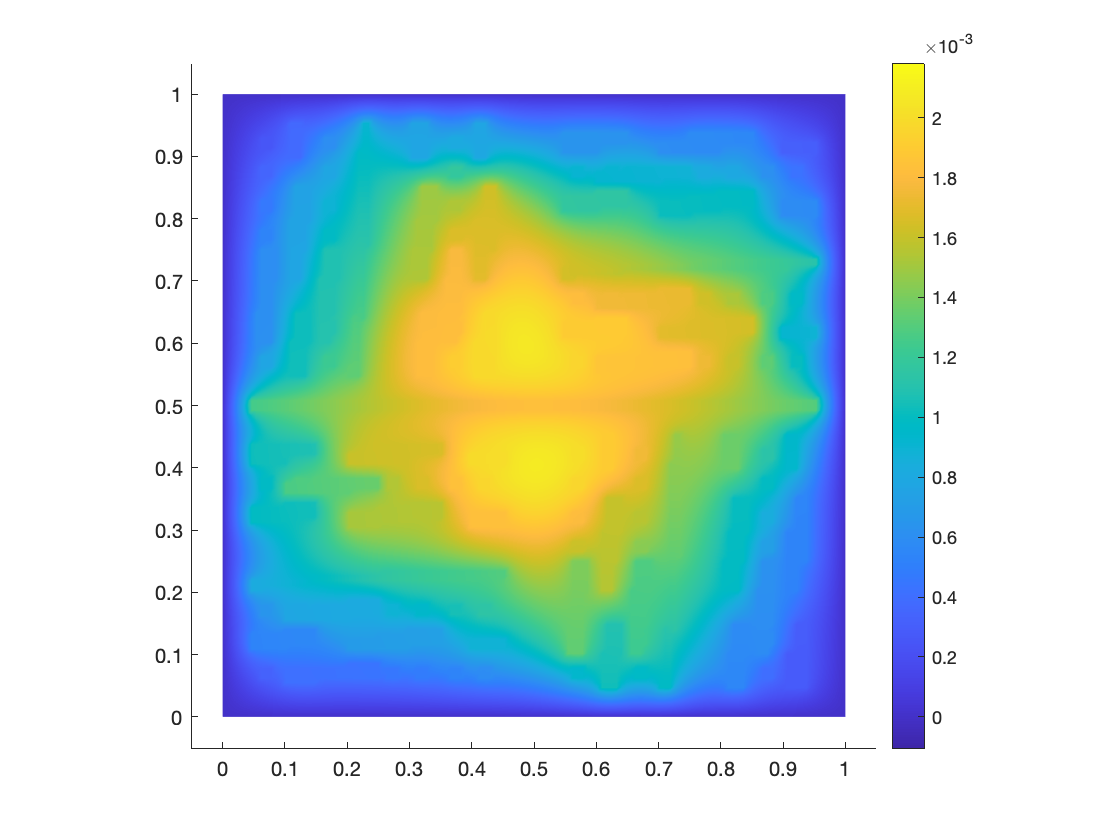}
  \caption{$p_{\tu{ms}}(T,\bfa{x})$ by CEM.}
  \label{sol_tdsingle_CEM}
\end{subfigure}
\hfill
   \begin{subfigure}{0.45\textwidth}
  \includegraphics[width=\textwidth]{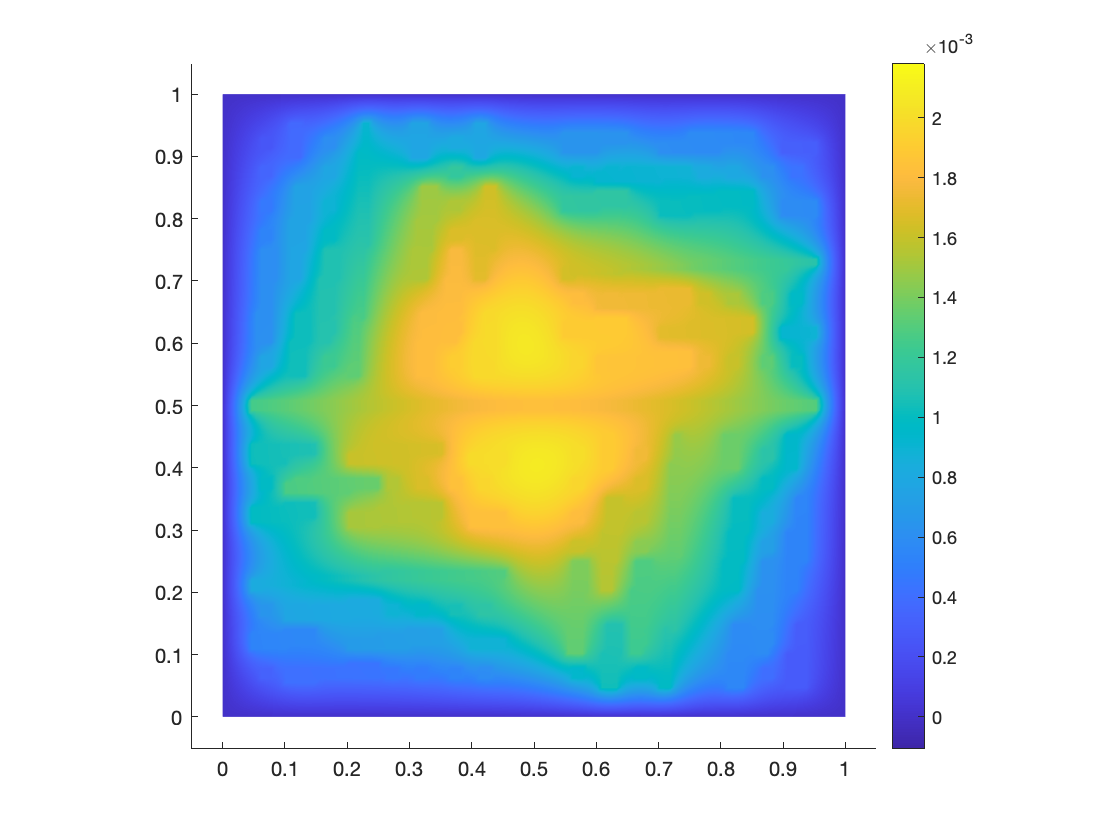}
  \caption{$p_h(T,\bfa{x})$ by FEM.}
  \label{sol_tdsingle_FEM}
 \end{subfigure}
 \caption{Problem (\ref{eq:original00}): Solutions $p(T,\bfa{x})$ obtained by CEM and FEM when $H=1/16\,.$}
 \label{fig:sol_tdsingle}
\end{figure}

\subsection{Experiments for dual-continuum Richards equations}
Benefiting from \cite{rtt21}, we consider the following steady-state problem of the form \eqref{eq:original} in the domain $\Omega\,:$ find $\bfa{p} = (p_1,p_2) \in 
[H_0^1(\Omega)]^2$ such that
\beq\label{nex}
\bsp
- \div \left(\frac{\kappa_1(\bfa{x})}{1+|p_1|} \nabla p_1 \right) + \frac{10}{1+|p_1|}(p_1-p_2) &= 1 \,,\\
- \div \left(\frac{\kappa_2(\bfa{x})}{1+|p_2|} \nabla p_2 \right) + \frac{10}{1+ |p_2|}(p_2-p_1) &= -1 \,,
\end{split}
\eeq
where it has zero Dirichlet boundary condition, and the configurations of high-contrast permeability fields $\kappa_1$ and $\kappa_2$ 
are shown in Figure \ref{fig:perm_steady}. The values in the yellow regions (channels) are higher than the values in the blue regions, and $\kappa_1(\bfa{x})\geq \kappa_2(\bfa{x})$ for all $\bfa{x}\in \Omega\,.$

As a special case of \eqref{eq:original} (with the given conditions there), the following problem is also considered.  That is, we investigate the Gardner-Basha model, which utilizes more intricate right-hand side functions and includes both sources and sinks 
(remark that the van Genuchten-Mualem model can work as well) \cite{15mualemv, santos2006hydraulic,rtt21}.  Employing the Gardner-Basha model in \cite{15mualemv}, we solely consider the unsaturated hydraulic conductivity's nonlinearity, and the volumetric water content function is assumed to be identity.  In the given domain $[0,T]\times \Omega\,,$ we seek solution $\bfa{p} = (p_1,p_2) \in 
[H_0^1(\Omega)]^2$ of the system
\beq\label{nex2}
\bsp
\frac{\partial p_1}{\partial t} - \div \left(\kappa_1(\bfa{x})K_r(p_1) \nabla p_1 \right)+ \frac{10^2}{1+|p_1|}(p_1-p_2) = f_1(\bfa{x}) \  \textrm{in} \ (0,T] \times \Omega \,,\\
\frac{\partial p_2}{\partial t} - \div \left(\kappa_2(\bfa{x}) K_r(p_2) \nabla p_2 \right)+ \frac{10^2}{1+ |p_2|}(p_2-p_1) = f_2(\bfa{x}) \  \textrm{in} \ (0,T] \times \Omega\,,
\end{split}
\eeq
with the Dirichlet boundary condition $p_i(t,\bfa{x})=0$ on $(0,T] \times \partial \Omega\,,$ and with the initial condition $p_i(0,\bfa{x})= 0$ in $\Omega\,.$  Here, the expression of relative hydraulic conductivity $K_{r}$ is
\begin{align}\label{Kr}
\begin{split}
	K_{r}(p) = e^{-\alpha_G|p|}\,,
\end{split}
\end{align}
where $\alpha_G$ is parameter characteristic of the soil pore size distribution.  The geometric mean of $\alpha_G$ is assumed to be $0.1\,.$  Fig.\ \ref{fig:perm_tdcoupled} describes the high-contrast $\kappa_i(\bfa{x})$ for $i = 1,2\,.$ We choose in the blue regions $\kappa_1(\bfa{x}) =10$ and $\kappa_2(\bfa{x}) =1\,,$ as well as in the yellow regions $\kappa_1(\bfa{x}) = 10^4$ and $\kappa_2(\bfa{x}) = 10\,.$  The specific source and sink functions are respectively provided by $f_1(\bfa{x}) = e^{x_1+x_2}$, $f_2(\bfa{x}) = -e^{x_1+x_2}\,.$

\begin{figure}[ht!]
\centering
\begin{subfigure}{0.45\textwidth}
  \includegraphics[width=\textwidth]{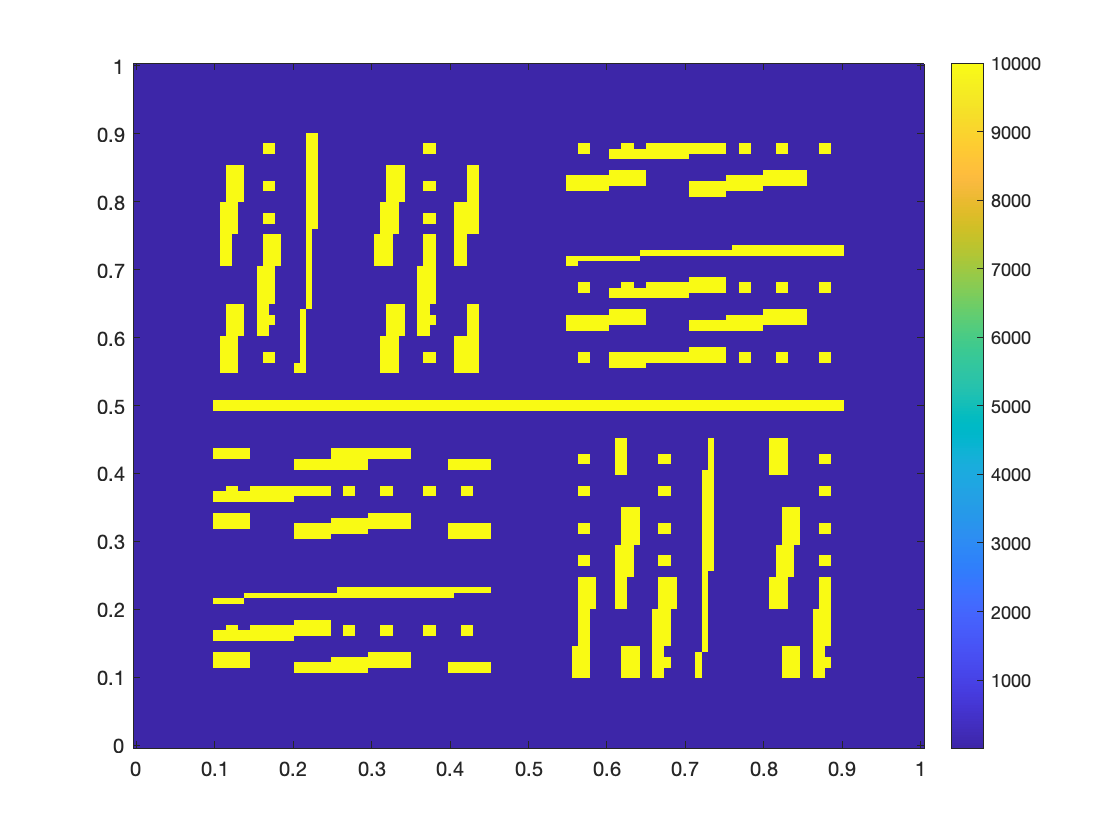}
  \caption{$\kappa_1(\bfa{x})\,.$}
  \label{a_1_steady}
\end{subfigure}
\hfill
   \begin{subfigure}{0.45\textwidth}
  \includegraphics[width=\textwidth]{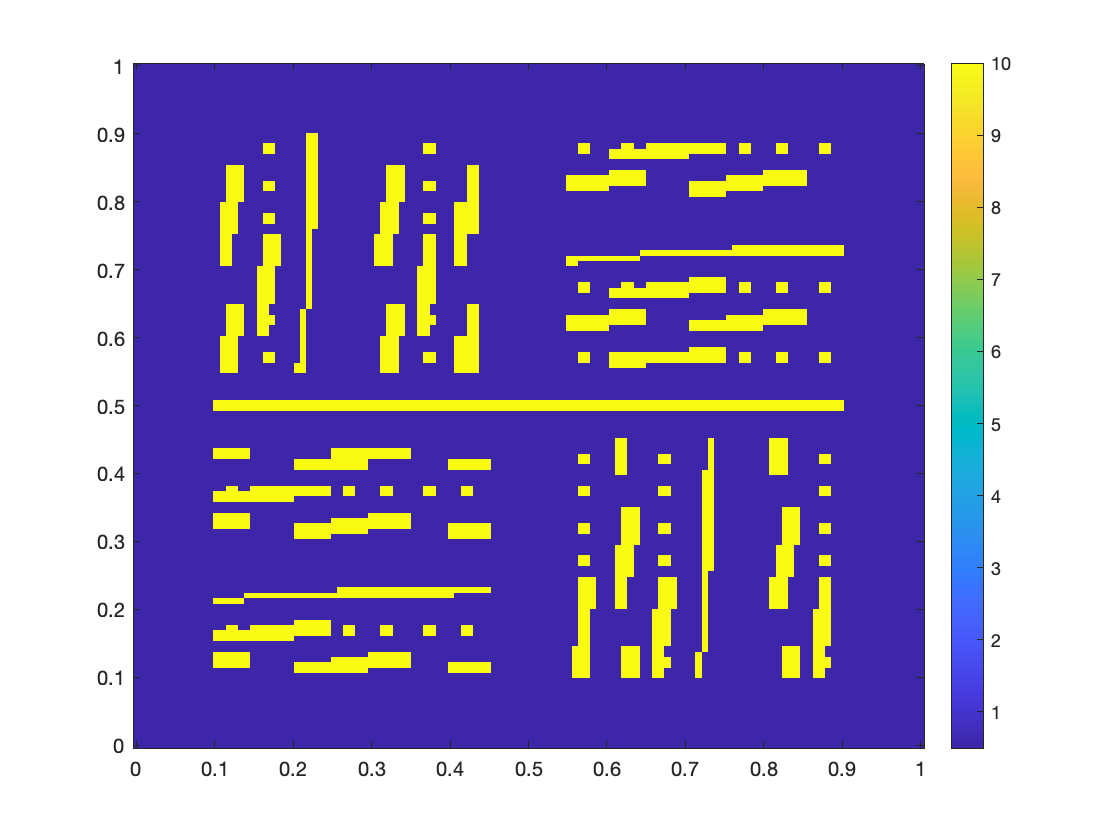}
  \caption{$\kappa_2(\bfa{x})\,.$}
  \label{a_2_steady}
 \end{subfigure}
 \caption{Problem (\ref{nex}): $\kappa_1(\bfa{x})=10^4$, $\kappa_2(\bfa{x})=10$ in the corresponding yellow regions (channels);  $\kappa_1(\bfa{x})=10$, $\kappa_2(\bfa{x})=0.5$ in the associated blue regions.}
 \label{fig:perm_steady}
\end{figure}
\begin{figure}[ht!]
\centering
\begin{subfigure}{0.45\textwidth}
  \includegraphics[width=\textwidth]{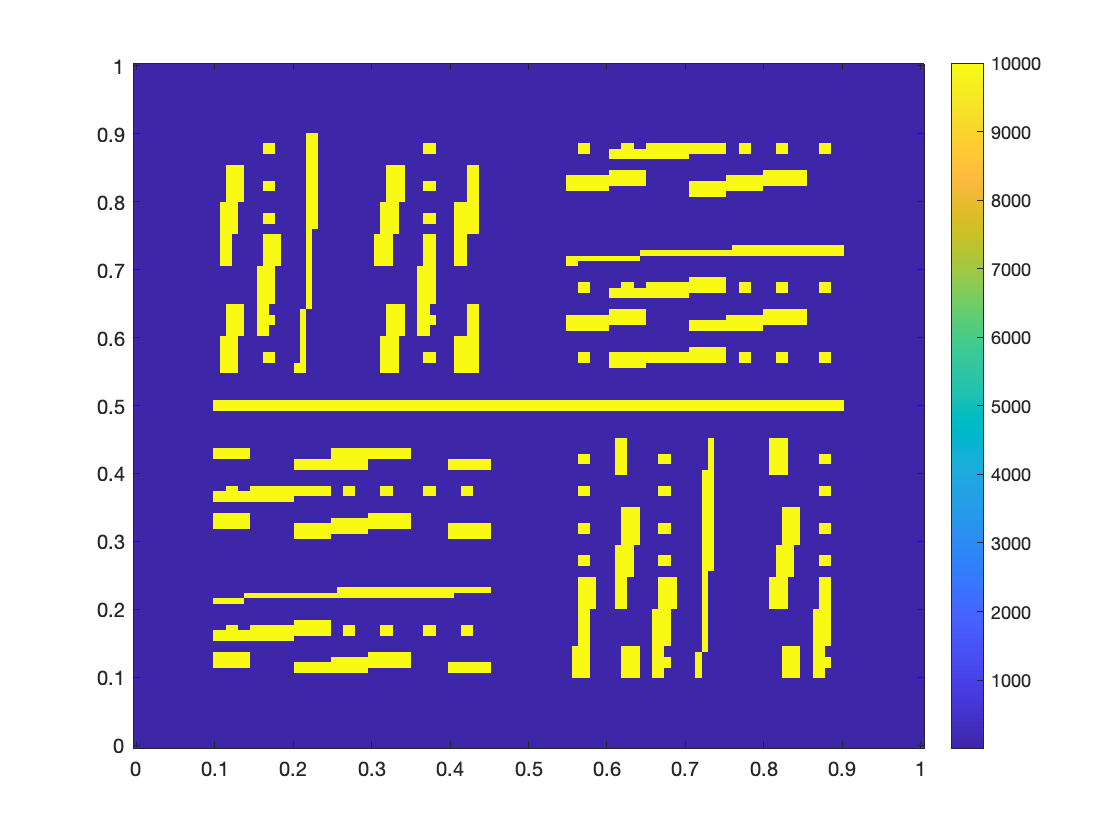}
  \caption{$\kappa_1(\bfa{x})\,.$}
  \label{a_1_tdcoupled}
\end{subfigure}
\hfill
   \begin{subfigure}{0.45\textwidth}
  \includegraphics[width=\textwidth]{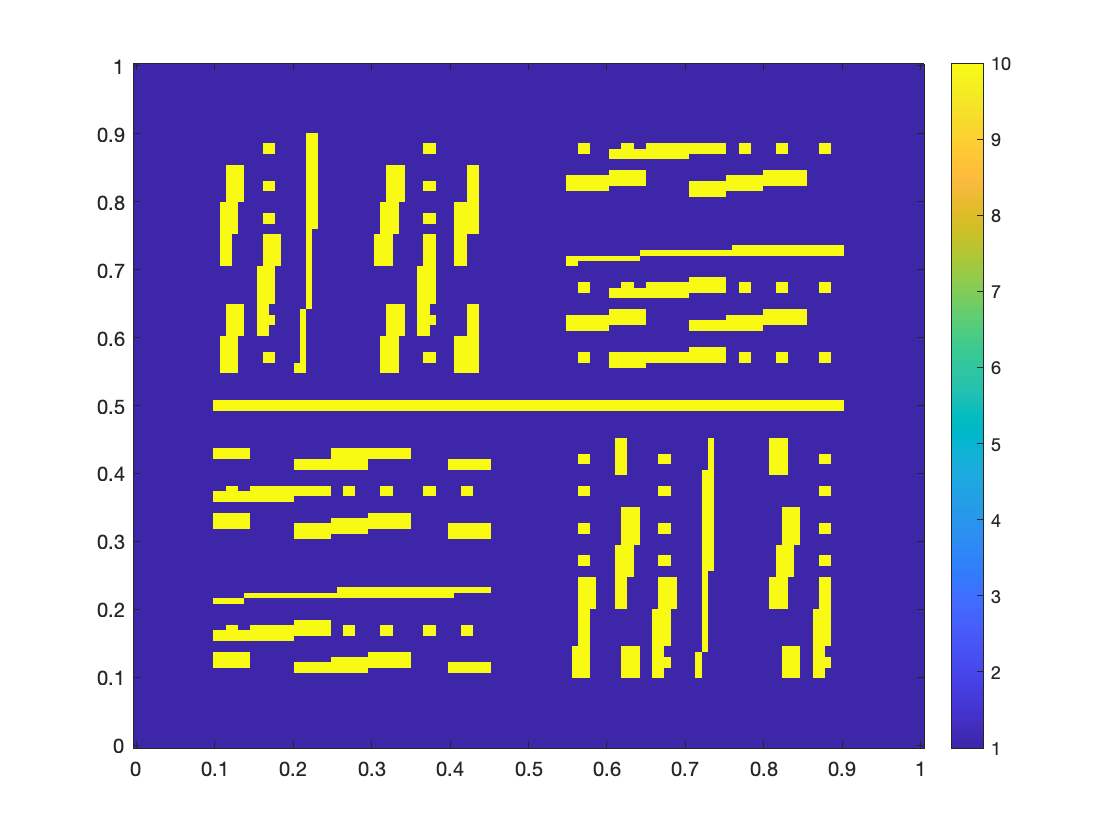}
  \caption{$\kappa_2(\bfa{x})\,.$}
  \label{a_2_tdcoupled}
 \end{subfigure}
 \caption{Problem (\ref{nex2}): $\kappa_1(\bfa{x})=10^4$, $\kappa_2(\bfa{x})=10$ in their yellow regions (channels);  $\kappa_1(\bfa{x})=10$, $\kappa_2(\bfa{x})=1$ in their blue regions.}
 \label{fig:perm_tdcoupled}
\end{figure}

We compute the numerical CEM solutions $\bfa{p}_{\tu{ms}}$ of (\ref{nex}) 
and (\ref{nex2}) (based on \eqref{ongmspicard}) as well as their respective FEM references $\bfa{p}_{\tu{h}}$ 
using \eqref{r1elh} without having to consider the time-step iterations for \eqref{nex}.  The following relative $L^2$ and $H^1$ errors are between the CEM solutions $\bfa{p}_{\tu{ms}}$ and their FEM references $\bfa{p}_{h}\,:$
\begin{align}\label{errors_dual}
	\begin{split}
		e^p_{L^2}=\frac{||\bfa{p}_{\tu{ms}}-\bfa{p}_{h}||_{\bfa{L}^2(\Omega)}}{||\bfa{p}_{h}||_{\bfa{L}^2(\Omega)}}, \
		e^p_{H^1}=\frac{||\nabla \bfa{p}_{\tu{ms}}-\nabla \bfa{p}_{h}||_{\mathbb{L}^2(\Omega)}}{||\nabla \bfa{p}_{h}||_{\mathbb{L}^2(\Omega)}}. 
	\end{split}
\end{align}
Table \ref{tab:steady} and \ref{tb:tdcoupled} present the relative $L^2$ and $H^1$ errors for problems (\ref{nex}) and (\ref{nex2}), respectively. We note that the total number of degrees of freedom of our multiscale method ($\textrm{dim}(V_{ms})$) depends on the coarse-grid size $H$ and the number of local multiscale basis functions. Those tables show clearly that the errors converge once the coarse-grid size $H$ is refined.  Also, according to the tables, increasing the number of local multiscale basis functions and oversampling layers help reduce the errors.  With $H = 1/4\,,$ the errors are relatively large especially for the time-dependent problem (\ref{nex2}), but they can be lower once more local multiscale basis functions are used. 
For relatively small coarse-grid size $H$, the error convergence tend to stagnate as we already have enough number of total degrees of freedom based on Table \ref{tab:steady} and \ref{tb:tdcoupled}.  
We observe that for small coarse-grid size, only few number of local multiscale basis functions are needed.  For the first continuum and at the final time $T=2\,,$ Figure \ref{fig:sol_tdcoupled} plots the solutions $p_1(T,\bfa{x})$ of (\ref{nex2}), obtained by the CEM and FEM when $H=1/16\,.$ 
Both solutions are almost identical throughout the entire domain.



\begin{table}[ht!]
\centering
\begin{subtable}{0.8\textwidth}
\centering
\begin{tabular}{|c|c|c||c||c|}
\hline
$H$ & $m$&$\textrm{dim}(V_{\tu{ms}})$ & $H^1$ error & $L^2$ error \\
\hline
$1/4$ & 3 &64& 5.6621\% & 1.2575\% \\
$1/8$ & 5 &256& 2.2713\%  & 0.2370\% \\
$1/16$ & 7 &1024&1.0382\% & 0.0736\%  \\
$1/32$ & 8 &4096& 0.4093\%  & 0.0135\%  \\
\hline
\end{tabular}
\caption{Errors with $4$ local basis functions for (\ref{nex}).}
\label{tab:steady_J4}
\end{subtable}
\begin{subtable}{0.8\textwidth}
\centering
\begin{tabular}{|c|c|c||c||c|}
\hline
$H$ & $m$&$\textrm{dim}(V_{\tu{ms}})$ & $H^1$ error & $L^2$ error \\
\hline
$1/4$ & 3 &80& 5.4959\% & 1.1861\% \\
$1/8$ & 5 &320& 2.0411\%  & 0.2042\% \\
$1/16$ & 7 &1280&1.0442\% & 0.0709\%  \\
$1/32$ & 8 &5120& 0.3896\%  & 0.0109\%  \\
\hline
\end{tabular}
\caption{Errors with $5$ local basis functions for (\ref{nex}).}
\label{tab:steady_J5}
\end{subtable}
\begin{subtable}{0.8\textwidth}
\centering
\begin{tabular}{|c|c|c||c||c|}
\hline
$H$ & $m$ &$\textrm{dim}(V_{\tu{ms}})$& $H^1$ error & $L^2$ error \\
\hline
$1/4$ & 3 &96& 4.6551\% & 0.9401\% \\
$1/8$ & 5 &384& 1.9344\%  & 0.1860\% \\
$1/16$ & 7 &1536&1.0086\% & 0.0647\%  \\
$1/32$ & 8 &6144& 0.3736\%  & 0.0092\%  \\
\hline
\end{tabular}
\caption{Errors with $6$ local basis functions for (\ref{nex}).}
\label{tab:steady_J6}
\end{subtable}
\vspace*{-3mm}
\caption{Relative $L^2$ and $H^1$ errors for (\ref{nex}) with different number of local basis functions; $\textrm{dim}(V_h) = 16129\,.$}
\label{tab:steady}
\end{table}

\begin{table}[ht!]
\centering
\begin{subtable}{0.8\textwidth}
\centering
\begin{tabular}{|c|c|c||c||c|}
\hline
$H$ & $m$ &$\textrm{dim}(V_{ms})$& $H^1$ error & $L^2$ error\\
\hline
$1/4$ & 3 &64& 15.7431\% & 5.7543\% \\
$1/8$ & 5 &256& 6.4665\% & 1.0124\% \\
$1/16$ & 7&1024 &1.7532\% & 0.1586\% \\
$1/32$ & 8 &4096& 0.6889\% &  0.0285\% \\
\hline
\end{tabular}
\caption{Errors with $4$ local basis functions for (\ref{nex2}).}
\label{tab:tdcoupled_J4}
\end{subtable}
\begin{subtable}{0.8\textwidth}
\centering
\begin{tabular}{|c|c|c||c||c|}
\hline
$H$ & $m$ &$\textrm{dim}(V_{ms})$& $H^1$ error & $L^2$ error\\
\hline
$1/4$ & 3 &80& 13.7012\% & 3.7477\% \\
$1/8$ & 5 &320& 4.5603\% & 0.6151\% \\
$1/16$ & 7 &1280&1.6636\% & 0.1481\% \\
$1/32$ & 8 &5120& 0.6465\% &  0.0223\% \\
\hline
\end{tabular}
\caption{Errors with $5$ local basis functions for (\ref{nex2}).}
\label{tab:tdcoupled_J5}
\end{subtable}
\begin{subtable}{0.8\textwidth}
\centering
\begin{tabular}{|c|c|c||c||c|}
\hline
$H$ & $m$ &$\textrm{dim}(V_{ms})$& $H^1$ error & $L^2$ error\\
\hline
$1/4$ & 3 &96& 10.7440\% & 2.5620\% \\
$1/8$ & 5 &384& 2.9378\% & 0.3349\% \\
$1/16$ & 7 &1536&1.5865\% & 0.1317\% \\
$1/32$ & 8 &6144& 0.6125\% &  0.0175\% \\
\hline
\end{tabular}
\caption{Errors with $6$ local basis functions for (\ref{nex2}).}
\label{tab:tdcoupled_J6}
\end{subtable}
\vspace*{-3mm}
\caption{Relative $L^2$ andf $H^1$ errors for (\ref{nex2}) with different number of local basis functions; $\textrm{dim}(V_h) = 16129\,.$}
\label{tb:tdcoupled}
\end{table}



\begin{figure}[ht!]
\centering
\begin{subfigure}{0.45\textwidth}
  \includegraphics[width=\textwidth]{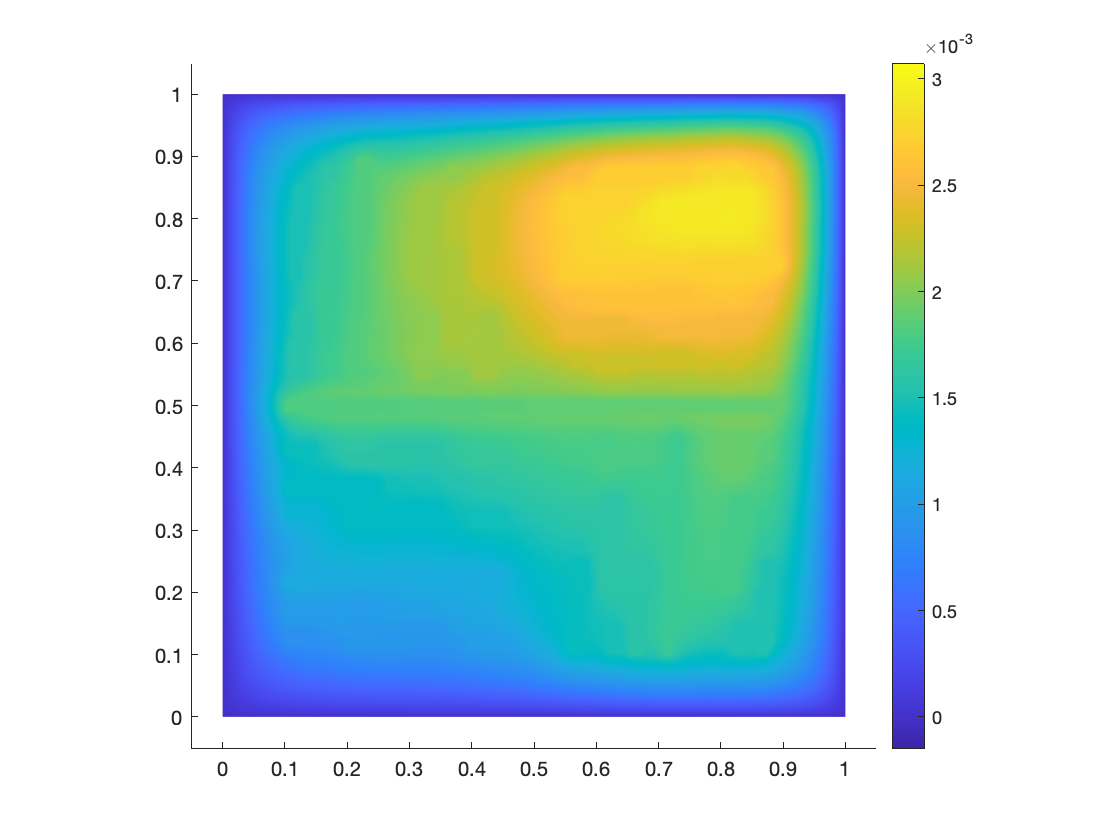}
  \caption{$p_{1,\tu{ms}}(T,\bfa{x})$ by CEM.}
  \label{P1_tdcoupled_CEM}
\end{subfigure}
\hfill
   \begin{subfigure}{0.45\textwidth}
  \includegraphics[width=\textwidth]{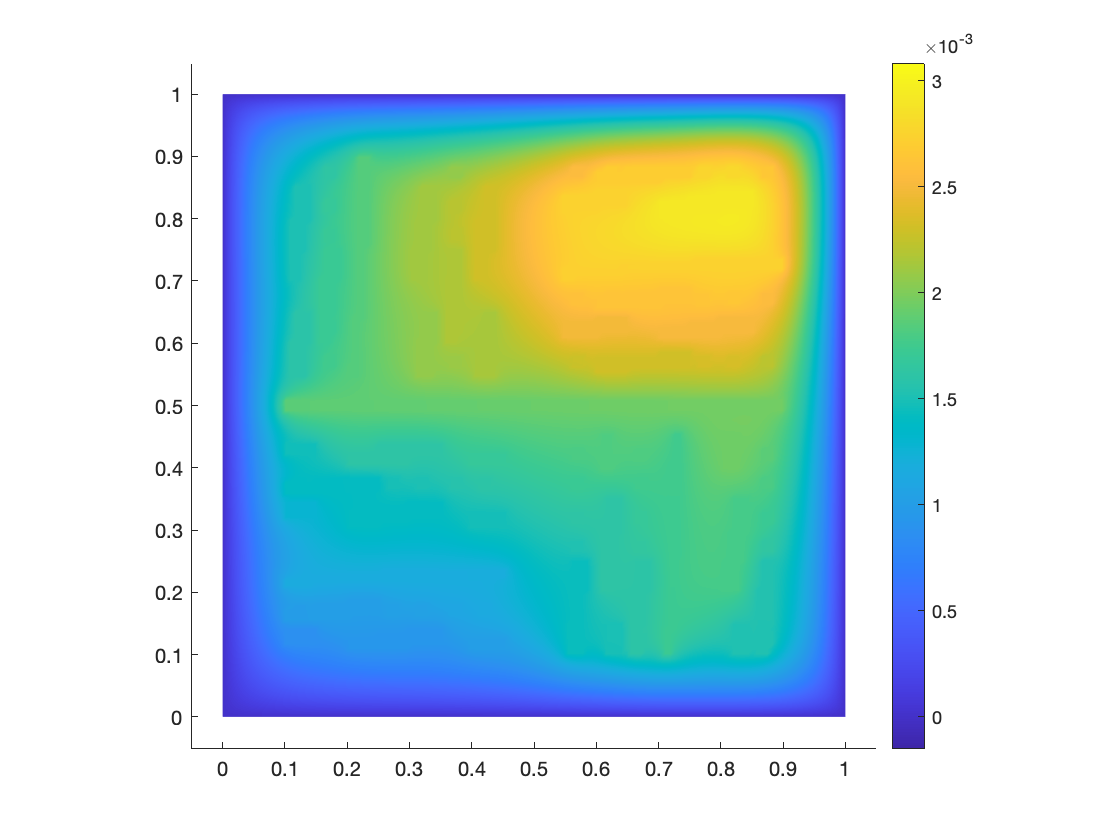}
  \caption{$p_{1,h}(T,\bfa{x})$ by FEM.}
  \label{P1_tdcoupled_FEM}
 \end{subfigure}
 \caption{Problem (\ref{nex2}): Solutions $p_1(T,\bfa{x})$ obtained by CEM and FEM when $T=2\,.$}
 \label{fig:sol_tdcoupled}
\end{figure}

\bigskip

\section{Conclusions}\label{sec:conclusions}
We present in this paper a methodology for handling issues from coupled system of multi-continuum nonlinear Richards equations, in complex fractured heterogeneous porous media, utilizing the constraint energy minimizing generalized multiscale finite element method (CEM-GMsFEM).  The basic concept is to discretize this system temporally via the implicit backward Euler method, then linearize it spatially by Picard iteration (with the required convergence indicator) at each time step until the ending time.  The CEM-GMsFEM is used in each Picard iteration to systematically create multiscale basis functions (with locally minimal energy) for pressure.  In order to do so, we propose two new sources of samples and solve local spectral problems via the GMsFEM to first build the local auxiliary multiscale basis functions, which are crucial for determining high-contrast channels. Second, employing the CEM through some constraints connected to the auxiliary functions, we solve an energy minimizing problem by oversampling technique, to establish localized multiscale basis functions. 
%
Our numerical results exhibits that the error converges with the coarse-grid size alone, and the method is very accurate.
Appendix \ref{cp} provides a theoretical proof for global convergence of the Picard iteration process.

\vspace{40pt}

\noindent \textbf{Acknowledgements.} 

Tina Mai's research was supported by RFBR and VAST under grant 21-51-54001 and by Duy Tan University under decision 5390/QD-DHDT.

This work was performed under the auspices of the U.S. Department of Energy by Lawrence Livermore National Laboratory under Contract DE-AC52-07NA27344 and LLNL-JRNL-833749.

\bigskip

\appendix 

\section{Global convergence of Picard linearization procedure}\label{cp}
%

We will prove the global convergence of Picard linearization process (given in Section \ref{pre}) for the system \eqref{r1elh} relied on \eqref{r1ed} and originated from \eqref{eq:original1}, with $i,l=1,2$ in this appendix, following \cite{rpicardc} (and thanks to J. Batista and A. Mazzucato, personal
communication, January 9, 2022).  The generalization to $i,l=1,\cdots,N$ is proved similarly.


In \eqref{eq:original} (and thus \eqref{r1e}, \eqref{r1ed} and \eqref{r1elh}), each hydraulic conductivity coefficient $\varkappa_i$ satisfies the assumption \eqref{Coercivity}, that is, $0 < \varkappa_1(\bfa{x}, p_1), \ \varkappa_2(\bfa{x}, p_2)\leq \overline{\varkappa}$ for some positive constant $\overline{\varkappa}\,.$  Each function $\varkappa_i$ is globally Lipschitz continuous with Lipschitz constant $L_{\varkappa_i}$ (without any dependence on $t$ and $\bfa{x}$).
Let $\bfa{p} = (p_1,p_2)\,,\bfa{p}_{il} = (p_i, p_l)\,,$ then the function $Q_{il}(\bfa{p}_{il}): = Q_{il}(\bfa{x},p_i,p_l)$ is nonlinear yet globally Lipschitz continuous with the Lipschitz constant $L_{Q_{il}}$ for every $i,l=1,2, i \neq l$ (without any dependence on $t$ and $\bfa{x}$).  Furthermore, we suppose that each $p_i$ is positive and that $Q_{il}$ is uniformly bounded above by some constant $\beta_{Q_{il}}$ ($\leq \overline{\beta}$ as from \eqref{Coercivity}).
The subscripts ($s+1$) and $h$ are eliminated from the Picard iteration \eqref{r1elh} for simplicity.

Then in $\bfa{V}_h\,,$ we obtain the following equality by subtracting \eqref{r1ed} from \eqref{r1elh} and picking correct $v_i = p_i^{n+1} - p_{i,s+1}\,:$
\begin{align}\label{in1}
\begin{split}
&\frac{1}{\tau} \|p_i^{n+1} - p_{i,s+1}\|^2
+ a_i(p_i^{n+1},p_i^{n+1} - p_{i,s+1};p_i^n) - a_i(p_{i,s+1},p_i^{n+1} - p_{i,s+1};p_{i,s+1}) \\
&=  - 
q_i(\bfa{p}^{n+1},p_i^{n+1} - p_{i,s+1};\bfa{p}^n) + 
q_i(\bfa{p}_{s+1},p_i^{n+1} - p_{i,s+1};\bfa{p}_{s+1})\,.
\end{split}
\end{align}
%

The left-hand side of \eqref{in1} is indicated by the following $LS$ using \eqref{ai}:
\begin{align}\label{lhs}
\begin{split}
LS &= \frac{1}{\tau} \|p_i^{n+1} - p_{i,s+1}\|^2 + (\varkappa(p_i^n) \nabla p_i^{n+1},\nabla(p_i^{n+1} - p_{i,s+1}))\\
& \qquad - (\varkappa(p_{i,s+1}) \nabla p_{i,s+1} , \nabla(p_i^{n+1} - p_{i,s+1}))\\
&= \frac{1}{\tau} \|p_i^{n+1} - p_{i,s+1}\|^2+ (\varkappa_i(p_i^n) \nabla (p_i^{n+1} - p_{i,s+1}),\nabla (p_i^{n+1} - p_{i,s+1})) \\
& \qquad + ((\varkappa_i(p_i^n) - \varkappa_i(p_{i,s+1}))\nabla p_{i,s+1} ,\nabla(p_i^{n+1} - p_{i,s+1}))\,.
\end{split}
\end{align}

The right-hand side of \eqref{in1} is represented by $RS$ employing \eqref{qi} as follows:

\vspace{-15pt}

\begin{align}\label{rhsin1}
\begin{split}
RS=&\sum_l \left\{ (-Q_{il}(\bfa{p}^n_{il})  \cdot (p_i^{n+1} - p_l^{n+1}), p_i^{n+1} - p_{i,s+1})\right.\\
&\left. \hspace{50pt}  + (Q_{il}(\bfa{p}_{il,s+1}) \cdot (p_{i,s+1} - p_{l,s+1}), p_i^{n+1} - p_{i,s+1})\right\}\\
& = \sum_l \left\{ (-Q_{il}(\bfa{p}_{il}^n) \cdot ((p_i^{n+1} - p_l^{n+1})-(p_{i,s+1} - p_{l,s+1})), p_i^{n+1} - p_{i,s+1})\right.\\
&\left. \hspace{40pt}  + ((-Q_{il}(\bfa{p}_{il}^n) + Q_{il}(\bfa{p}_{il,s+1}))\cdot (p_{i,s+1} - p_{l,s+1}), p_i^{n+1} - p_{i,s+1})\right\}\,.
\end{split}
\end{align}

For $i,l=1,2\,, i\neq l \,, $ one notes that
\begin{equation}\label{nvvf}
\|p^n_i-p^n_l\| \leq \|(p^n_i,p^n_l)\| = \|\bfa{p}^n\| \,, \qquad \|p_{l}^{n} - p_{l,s+1}\| \leq \| \bfa{p}^{n} - \bfa{p}_{s+1}\| \leq
\displaystyle \sum_{l} \|p_{l}^{n} - p_{l,s+1}\|\,.
\end{equation}
Assume there are respectively sufficiently large and small constants $D_l > 0$ and $M_l > 0$ such that $M_l \leq \|p_l^{n+1} - p_{l,s+1}\| \leq \hat{C}  \|\nabla (p_l^{n+1} - p_{l,s+1})\| \leq D_l\,,$ (for $l=1,2$), where the constant $\hat{C}$ depends only on $\Omega\,.$
Also, assume that each $p_{i,s+1}(t, \cdot) \in C_c^{\infty}(\Omega)$ so that $\|  \bfa{p}_{s+1} \|_{\infty} \leq \hat{D} \|\nabla  \bfa{p}_{s+1} \|_{\infty}\,,$ where $\hat{D}$ is the distance between the two parallel hyperplanes bounding $\Omega\,.$  
Let $\bfa{U}$ be the exact solution of the problem at hand \eqref{eq:original}.  Then, we obtain the following inequalities by employing the error estimate in \cite{GalerkinFEM} (Theorem 1.5):
\begin{align}\label{re5}
\begin{split}
\frac{1}{\hat{D}} \|  \bfa{p}_{s+1} \|_{\infty} \leq \| \nabla \bfa{p}_{s+1} \|_{\infty} \leq  \|\nabla(\bfa{U}(t_{s+1}) - \bfa{p}_{s+1}) \|_{\infty} + \|\nabla \bfa{U}(t_{s+1})\|_{\infty} & \leq  \overline{C}(\bfa{U})(h + \tau) + \|\nabla \bfa{U}(t_{s+1})\|_{\infty} = \overline{M}\,,
\end{split}
\end{align}
for some constant $\overline{C}(\bfa{U})$ depending on $\bfa{U}\,.$ 

We therefore get from \eqref{in1}, \eqref{lhs}, \eqref{rhsin1}, 
Young's inequality, and \eqref{nvvf} that for $i,l=1,2, i\neq l\,,$
\begin{align}\label{in1xi}
\begin{split}
&\frac{1}{\tau} \|p_i^{n+1} - p_{i,s+1}\|^2 + \underline{\varkappa}  \|\nabla (p_i^{n+1} - p_{i,s+1})\|^2 
\\
& \quad \leq L_{\varkappa_i} \|p_i^n - p_{i,s+1}\| \, \|\nabla p_{i,s+1} \|_{\infty} \|\nabla(p_i^{n+1} - p_{i,s+1})\|\\
& \qquad  + \sum_l \left\{ \beta_{Q_{il}} \|(p_i^{n+1} - p_l^{n+1})-(p_{i,s+1} - p_{l,s+1})\| \ \|p_i^{n+1} - p_{i,s+1}\|\right.\\
&\left. \hspace{55pt}  + \|(-Q_{il}(\bfa{p}_{il}^n) + Q_{il}(\bfa{p}_{il,s+1})\| \ \|p_{i,s+1} - p_{l,s+1}\| \ \|p_i^{n+1} - p_{i,s+1}\|\right\}\\
& \quad
\leq \frac{L^2_{\varkappa_i} \|\nabla p_{i,s+1} \|_{\infty}^2}{2 \underline{\varkappa}} \|p_i^n - p_{i,s+1}\|^2 + \frac{\underline{\varkappa}}{2} \|\nabla(p_i^{n+1} - p_{i,s+1})\|^2\\
& \qquad + \beta_{Q_{il}} \, \| \bfa{p}^{n+1} - \bfa{p}_{s+1}\| \, \|p_i^{n+1} - p_{i,s+1}\|\\
& \qquad + |\Omega| L_{Q_{il}}\, \| \bfa{p}^n - \bfa{p}_{s+1}\|  \, (\|\bfa{p}_{s+1}\|_{\infty} )
\|p_i^{n+1} - p_{i,s+1}\|\,.
\end{split}
\end{align}
Equivalently,
\begin{align}\label{in1xi2}
\begin{split}
& \left(\frac{1}{\tau} + \frac{\underline{\varkappa}}{2\hat{C}} \right) \|p_i^{n+1} - p_{i,s+1}\| \\
& \quad \leq  
\frac{L^2_{\varkappa_i} \|\nabla p_{i,s+1} \|_{\infty}^2}{2 \underline{\varkappa}} \|p_i^{n} - p_{i,s+1}\| \\
& \qquad + 
|\Omega| L_{Q_{il}}\, \| \bfa{p}^n - \bfa{p}_{s+1}\|  \, (\|\bfa{p}_{s+1}\|_{\infty} )
+ \beta_{Q_{il}} \, \| \bfa{p}^{n+1} - \bfa{p}_{s+1}\| \,.
\end{split}
\end{align}

Now, let $D = \max\{D_l\}\,,$ $L= \max \{L_{\varkappa_i}, |\Omega| L_{Q_{il}}\}\,,$ and recall that $\overline{\beta} = \max \{\beta_{Q_{il}}\}\,.$  
Letting $i=1,2$ in \eqref{in1xi2}, then summing up the resulting inequalities and benefiting from \eqref{nvvf}, we have 
\begin{align}\label{in1xi3}
\begin{split}
&\left(\frac{1}{\tau} + \frac{\underline{\varkappa}}{2\hat{C}} \right) \|\bfa{p}^{n+1} - \bfa{p}_{s+1}\| \\
&\quad \leq \frac{L^2 \|\nabla \bfa{p}_{s+1} \|_{\infty}^2}{ \underline{\varkappa}}   \, \| \bfa{p}^n - \bfa{p}_{s+1}\| \\
& \qquad +   |\Omega|(L_{Q_{12}}+L_{Q_{21}}) \, \| \bfa{p}^n - \bfa{p}_{s+1}\|  \, \|\bfa{p}_{s+1}\|_{\infty} 
+ (\beta_{Q_{12}} + \beta_{Q_{21}}) \, \| \bfa{p}^{n+1} - \bfa{p}_{s+1}\| \\
& \quad \leq  \left(\frac{L^2}{ \underline{\varkappa}} \|\nabla \bfa{p}_{s+1} \|_{\infty}^2  + 2L\hat{D} \|\nabla \bfa{p}_{s+1}\|_{\infty} \right)\| \bfa{p}^n - \bfa{p}_{s+1}\|   + 2\overline{\beta} \, \| \bfa{p}^{n+1} - \bfa{p}_{s+1}\|\,.
\end{split}
\end{align}
Here, we assume that $\tau$ is very small such that the left-hand side of \eqref{in1xi3} is much larger than the second term of the last right-hand side of \eqref{in1xi3}.  After rearranging the later inequality of \eqref{in1xi3} as well as applying \eqref{re5}, we reach
\begin{align}\label{re3}
\begin{split}
\|\bfa{p}^{n+1} - \bfa{p}_{s+1}\| 
\leq L \frac{2\tau \hat{C}}{ \underline{\varkappa}  (2\hat{C}+ \tau \underline{\varkappa} - 4 \overline{\beta} \tau \hat{C})} (L\overline{M}^2 + 2\underline{\varkappa} \hat{D} \overline{M})  \| \bfa{p}^n - \bfa{p}_{s+1}\|\,.
\end{split}
\end{align}
%
%
%

Upon redefining constants appropriately in \eqref{re3}, we eventually get
\begin{equation}\label{re6}
\| \bfa{p}^{n+1} - \bfa{p}_{s+1}\| \leq 
\frac{C L} {\underline{\varkappa}} \ \frac{\tau }{1 + \dfrac{( \underline{\varkappa} - 4 \overline{\beta} \hat{C})\tau}{2\hat{C}}} (\tilde{C}(h+\tau) +1)^2 \| \bfa{p}^{n} - \bfa{p}_{s+1}\|: = \lambda \| \bfa{p}^n - \bfa{p}_{s+1}\|\,,
\end{equation}
for some positive constants $C, \tilde{C}\,.$
With sufficiently small $\tau$ and $h\,,$ the coefficient $\lambda$ will be less than $1\,,$ implying that the algorithm converges.  Specifically, $\lambda \to 0$ when $h \to 0$ and $\tau \to 0$ at the same time.


\bibliographystyle{plain}
\bibliography{r1,r2}

\end{document}